\documentclass[a4paper,12pt]{article}

\usepackage{amsmath,amssymb,amsthm,amscd}

\newcommand{\Fg}{\mathfrak{g}}
\newcommand{\Fh}{\mathfrak{h}}

\newcommand{\BB}{\mathbb{B}}
\newcommand{\BP}{\mathbb{P}}
\newcommand{\BZ}{\mathbb{Z}}
\newcommand{\BQ}{\mathbb{Q}}
\newcommand{\BR}{\mathbb{R}}

\newcommand{\CB}{\mathcal{B}}
\newcommand{\CS}{\mathcal{S}}
\newcommand{\CT}{\mathcal{T}}

\newcommand{\wt}{\mathop{\rm wt}\nolimits}
\newcommand{\Ht}{\mathop{\rm ht}\nolimits}
\newcommand{\ch}{\mathop{\rm ch}\nolimits}
\newcommand{\cl}{\mathop{\rm cl}\nolimits}
\newcommand{\Wt}{\mathop{\rm Wt}\nolimits}
\newcommand{\Cat}{\mathop{\rm Cat}\nolimits}

\newcommand{\vpi}{\varpi}
\newcommand{\ve}{\varepsilon}
\newcommand{\vp}{\varphi}

\newcommand{\ti}[1]{\widetilde{#1}}
\newcommand{\ud}[1]{\underline{#1}}

\makeatletter
\@addtoreset{equation}{subsection}

\renewcommand\section{\@startsection{section}{1}{0pt}
{-3.5ex plus -1ex minus -.2ex}{1.0ex plus .2ex}{\large\bf}}
\renewcommand\subsection{\@startsection{subsection}{1}{0pt}
{2.5ex plus 1ex minus .2ex}{-1em}{\bf}}

\theoremstyle{plain}
\newtheorem{thm}{Theorem}[subsection]
\newtheorem{lem}[thm]{Lemma}
\newtheorem{prop}[thm]{Proposition}
\newtheorem{cor}[thm]{Corollary}
\newtheorem{ithm}{Theorem}
\newtheorem*{icor}{Corollary}

\theoremstyle{definition}

\theoremstyle{remark}
\newtheorem{rem}[thm]{Remark\rm}
\newtheorem*{irem}{Remark}

\begin{document}

\setlength{\baselineskip}{18pt}
\setcounter{section}{-1}


\title{{\Large\bf Path Model for a Level Zero 
Extremal Weight Module over a Quantum Affine Algebra%
\footnote{To be revised.}}}
\author{Satoshi Naito \\ 
\small Institute of Mathematics, University of Tsukuba, \\
\small Tsukuba, Ibaraki 305-8571, Japan \ 
(e-mail: {\tt naito@math.tsukuba.ac.jp})
\\[2mm] and \\[2mm]
Daisuke Sagaki \\ 
\small Institute of Mathematics, University of Tsukuba, \\
\small Tsukuba, Ibaraki 305-8571, Japan \ 
(e-mail: {\tt sagaki@math.tsukuba.ac.jp})}
\date{}
\maketitle

\begin{abstract}
We give a path model for a level zero extremal weight module 
over a quantum affine algebra. By using this result, we prove 
a branching rule for an extremal weight module with respect 
to a Levi subalgebra. Furthermore, we also show a 
decomposition rule of Littelmann type 
for the concatenation of path models for an 
integrable highest weight module and 
a level zero extremal weight module in the case where 
the extremal weight is minuscule.
\end{abstract}

%
\section{Introduction.}
\label{sec:intro}

Let $\Fg$ be a symmetrizable Kac--Moody algebra
with $\Fh$ the Cartan subalgebra and $W$ the Weyl group. 
We fix an integral weight 
lattice $P \subset \Fh^{\ast}$ that contains 
all simple roots $\Pi$ of $\Fg$. 
Let $\lambda \in P$ be an integral weight. 
In \cite{L1} and \cite{L2}, Littelmann 
introduced the notion of Lakshmibai--Seshadri paths 
of shape $\lambda$, which are piecewise linear, continuous 
maps $\pi:[0,1] \rightarrow P$ parametrized by pairs 
of a sequence of elements of $W\lambda$ and 
a sequence of rational numbers satisfying a certain condition, 
called the chain condition. Denote by $\BB(\lambda)$ 
the set of Lakshmibai--Seshadri paths of shape $\lambda$. 
Littelmann proved that $\BB(\lambda)$ has a normal 
crystal structure in the sense of \cite{Kas3}, 
and that if $\lambda$ is a dominant 
integral weight, then the formal sum
$\sum_{\pi \in \BB(\lambda)}e(\pi(1))$ is equal to the character 
$\ch L(\lambda)$ of 
the integrable highest weight $\Fg$-module $L(\lambda)$ 
of highest weight $\lambda$. Then he conjectured that 
$\BB(\lambda)$ for dominant $\lambda \in P$ would be 
isomorphic to the crystal base of the integrable highest weight 
module of highest weight $\lambda$ as crystals. 
This conjecture was affirmatively proved independently by 
Kashiwara \cite{Kas4} and Joseph \cite{J}. 

In \cite{Kas2} and \cite{Kas5}, Kashiwara 
introduced an extremal weight module $V(\lambda)$
of extremal weight $\lambda \in P$ over a quantized universal 
enveloping algebra $U_{q}(\Fg)$, and showed that it has a 
crystal base $\CB(\lambda)$. The extremal weight module is 
a natural generalization of an integrable 
highest (lowest) weight module. In fact, 
we know from \cite[\S8]{Kas2} 
that if $\lambda \in P$ is dominant (resp. anti-dominant), 
then the extremal weight module $V(\lambda)$ is 
isomorphic to the integrable highest (resp. 
lowest) weight module of highest (resp. lowest) weight $\lambda$, 
and the crystal base $\CB(\lambda)$ of $V(\lambda)$ 
is isomorphic to the crystal base of the integrable highest (resp. 
lowest) weight module as crystals. 

Now, assume that $\Fg$ is of affine type. Let $I$ be the index set 
of the simple roots $\Pi$ of $\Fg$, and fix a special vertex $0 \in I$
as in \cite[\S5.2]{Kas5}. In this paper, as an extension of the 
isomorphism theorem due to Kashiwara and Joseph, 
we prove that if $\lambda$ 
is a fundamental basic weight $\vpi_{i} \in P$ for 
$i \in I_{0}:=I \setminus \{0\}$ (see \cite[\S5.2]{Kas5}; 
note that $\vpi_{i}$ is not dominant), 
then the connected component\footnote{The authors have checked that 
the crystal graph of $\BB(\vpi_{i})$ is connected for 
some special $\vpi_{i}$'s. For example, 
if $\vpi_{i}$ is minuscule, then $\BB(\vpi_{i})$ 
is connected. The connectedness 
of the crystal $\BB(\vpi_{i})$ will 
be studied in a forthcoming paper.} $\BB_{0}(\vpi_{i})$
of $\BB(\vpi_{i})$ containing $\pi_{\vpi_{i}}(t):=t\vpi_{i}$ is 
isomorphic to the crystal base $\CB(\vpi_{i})$ of the extremal weight 
module $V(\vpi_{i})$ as crystals. Namely, we prove the following:

\begin{ithm}
Assume that $\Fg$ is of affine type. 
There exists a unique isomorphism $\Phi_{\vpi_{i}}:
\CB(\vpi_{i}) \stackrel{\sim}{\rightarrow} 
\BB_{0}(\vpi_{i})$ of crystals such that $\Phi_{\vpi_{i}}
(u_{\vpi_{i}})=\pi_{\vpi_{i}}$, where $u_{\vpi_{i}} \in 
\CB(\vpi_{i})$ is the unique extremal weight element of 
extremal weight $\vpi_{i}$. 
\end{ithm}

In \cite{Kas5}, Kashiwara conjectured that if $\lambda=
\sum_{i \in I_{0}} m_{i}\vpi_{i}$, then there exists a natural 
embedding of crystals from the connected component 
$\CB_{0}(\lambda)$ of $\CB(\lambda)$ containing 
the extremal weight element $u_{\lambda}$ of weight $\lambda$ 
into the tensor product $\bigotimes_{i \in I_{0}} 
\CB(\vpi_{i})^{\otimes m_{i}}$ that maps $u_{\lambda}$
to $\bigotimes_{i \in I_{0}} u_{\vpi_{i}}^{\otimes m_{i}}$. 
Beck and Nakajima proved in \cite{B}, \cite{BN}, and \cite{N} 
that this conjecture is true. Because the tensor product of 
crystals corresponds to the concatenation of path models 
(cf. \cite[\S1 and Lemma 2.7]{L2}), 
we get the following corollary by combining 
the above result of Beck and Nakajima and Theorem 1. 

\begin{icor}
Let $\lambda=\sum_{i \in I_{0}} m_{i}\vpi_{i} \in P$ 
with $m_{i} \in \BZ_{\ge 0}$, and let
$\BB_{0}(\ud{\lambda})$ be the connected component of 
the concatenation
$\Cat_{i \in I_{0}} \BB_{0}(\vpi_{i})^{\ast m_{i}}$
of path models $\BB_{0}(\vpi_{i})$'s containing the path 
$\pi_{\ud{\lambda}}:=\Cat_{i \in I_{0}}
\pi_{\vpi_{i}}^{\ast m_{i}}$, i.e., that of 
$\Cat_{i \in I_{0}} \BB(\vpi_{i})^{\ast m_{i}}$.
Then there exists a unique isomorphism $\Phi_{\ud{\lambda}}:
\CB_{0}(\lambda) \stackrel{\sim}{\rightarrow} 
\BB_{0}(\ud{\lambda})$ of crystals 
that maps $u_{\lambda}$ to $\pi_{\ud{\lambda}}$. 
\end{icor}

Let $\Fg_{S}$ be a Levi subalgebra corresponding to 
a proper subset $S$ of the index set $I$. 
By restriction, we can regard the crystals
$\BB(\vpi_{i})$ and $\BB_{0}(\vpi_{i})$ for $\Fg$ as crystals 
for $\Fg_{S}$. We show the following branching rule for
$\BB(\vpi_{i})$ and $\BB_{0}(\vpi_{i})$ as crystals for $\Fg_{S}$:

\begin{ithm}
As crystals for $\Fg_{S}$, $\BB(\vpi_{i})$ and $\BB_{0}(\vpi_{i})$ 
decompose as follows\,{\rm:}
\begin{equation*}
\BB(\vpi_{i}) \cong \bigsqcup_{
\begin{subarray}{c}
\pi \in \BB(\vpi_{i}) \\
\pi\text{\rm: $\Fg_{S}$-dominant}
\end{subarray}}
\BB_{S}(\pi(1)), \qquad 
\BB_{0}(\vpi_{i}) \cong \bigsqcup_{
\begin{subarray}{c}
\pi \in \BB_{0}(\vpi_{i}) \\
\pi\text{\rm: $\Fg_{S}$-dominant}
\end{subarray}}
\BB_{S}(\pi(1)).
\end{equation*}
where $\BB_{S}(\lambda)$ is the set of Lakshmibai--Seshadri paths
of shape $\lambda$ for $\Fg_{S}$, and $\pi \in \BB(\vpi_{i})$
is said to be $\Fg_{S}$-dominant if $(\pi(t))(\alpha^{\vee}_{i}) \ge 0$
for all $t \in [0,1]$ and $i \in S$.
\end{ithm}

We also show that the extremal weight module $V(\vpi_{i})$ 
of extremal weight $\vpi_{i}$ is completely reducible as a 
$\Fg_{S}$-module. As an application of Theorems 1 and 2 above, 
we obtain the following branching rule for $V(\vpi_{i})$: 

\begin{ithm}
The extremal weight module $V(\vpi_{i})$ of extremal weight $\vpi_{i}$ 
is completely reducible as a $U_{q}(\Fg_{S})$-module, and 
the decomposition of $V(\vpi_{i})$ as a $U_{q}(\Fg_{S})$-module 
is given by\,{\rm:}
\begin{equation*}
V(\vpi_{i}) \cong \bigoplus_{
\begin{subarray}{c}
\pi \in \BB_{0}(\vpi_{i}) \\
\pi\text{\rm: $\Fg_{S}$-dominant}
\end{subarray}}
V_{S}(\pi(1)), 
\end{equation*}
where $V_{S}(\lambda)$ is the integrable highest weight module 
of highest weight $\lambda$ over the quantized universal 
enveloping algebra $U_{q}(\Fg_{S})$ of the Levi subalgebra $\Fg_{S}$. 
\end{ithm}

Assume that $\vpi_{i}$ is minuscule, i.e., 
$\vpi_{i}(\alpha^{\vee}) \in \bigl\{\pm1,0\bigr\}$
for every dual real root $\alpha^{\vee}$ of $\Fg$.
Then we can check that $\BB(\vpi_{i})$ is connected, and hence 
$\BB(\vpi_{i})=\BB_{0}(\vpi_{i})$. 
In this case, we get the following decomposition rule of 
Littelmann type for the concatenation
$\BB(\lambda) \ast \BB(\vpi_{i})$. Here we note 
that unlike Theorems 2 and 3, 
this theorem does not necessarily 
imply the decomposition rule for tensor products of
corresponding $U_{q}(\Fg)$-modules.

\begin{ithm}
Let $\lambda$ be a dominant integral weight which is not a multiple
of the null root $\delta$ of $\Fg$, and assume that 
$\vpi_{i}$ is minuscule. Then 
the concatenation $\BB(\lambda) \ast \BB(\vpi_{i})$ 
decomposes as follows\,{\rm:}
\begin{equation*}
\BB(\lambda) \ast \BB(\vpi_{i})
 \cong {\displaystyle\bigsqcup_{
\begin{subarray}{c}
\pi \in \BB(\vpi_{i}) \\
\pi\text{\rm: $\lambda$-dominant}
\end{subarray}}}
\BB(\lambda+\pi(1)), 
\end{equation*}
where $\pi \in \BB(\vpi_{i})$ is said to be $\lambda$-dominant if 
$(\lambda+\pi(t))(\alpha^{\vee}_{i}) \ge 0$
for all $t \in [0,1]$ and $i \in I$.
\end{ithm}

\begin{irem}
The reader should compare Theorems 1 and 4 with 
the corresponding results \cite[Theorems 1.5 and 1.6]{G}
of Greenstein for bounded modules.
\end{irem}

This paper is organized as follows. In \S\ref{sec:prenot}, 
we fix our notation, and recall some basic facts about crystal 
bases. In \S\ref{sec:pre01} and \S\ref{sec:pre02}, we introduce 
some tools for crystal bases and path models, 
which will be needed in the proof of our isomorphism 
theorem (Theorem 1). In \S\ref{sec:main}, 
we will prove our main results mentioned above. 

%
\section{Preliminaries and Notation.}
\label{sec:prenot}

%
\subsection{Quantized universal enveloping algebras.}
\label{subsec:quea}

Let $A=(a_{ij})_{i,j \in I}$ be a symmetrizable 
generalized Cartan matrix, and $\Fg:=\Fg(A)$ the 
Kac--Moody algebra over $\BQ$ associated to 
the generalized Cartan matrix $A$. Denote by 
$\Fh$ the Cartan subalgebra, 
by $\Pi:=\bigl\{\alpha_{i}\bigr\}_{i \in I} 
\subset \Fh^{\ast}$ and $\Pi^{\vee}:=
\bigl\{\alpha^{\vee}_{i}\bigr\}_{i \in I} \subset \Fh$ 
the set of simple roots and simple coroots, 
and by $W=\langle r_{i} \mid 
i \in I \rangle$ the Weyl group. We take (and fix) 
an integral weight lattice $P \subset \Fh^{\ast}$ 
such that $\alpha_{i} \in P$ for all $i \in I$. 

Denote by $U_{q}(\Fg)$ the quantized universal enveloping 
algebra of $\Fg$ over the field $\BQ(q)$ of rational 
functions in $q$, and 
by $U_{q}^{-}(\Fg)$ (resp. $U_{q}^{+}(\Fg)$) the negative 
(resp. positive) part of $U_{q}(\Fg)$. We denote by 
$\ti{U}_{q}(\Fg)=\bigoplus_{\lambda \in P} 
U_{q}(\Fg)a_{\lambda}$ the modified quantized universal 
enveloping algebra of $\Fg$, where $a_{\lambda}$ is 
a formal element of weight $\lambda$ (cf. \cite[\S1.2]{Kas2}).

%
\subsection{Affine Lie algebras.}
\label{subsec:affine}
In some parts of this paper, we will assume that $\Fg$ is of 
affine type. In this case, we use the following notation. 
Let 
\begin{equation}
\delta=\sum_{i \in I} a_{i}\alpha_{i} \in \Fh^{\ast}
\qquad \text{ and } \qquad 
c=\sum_{i \in I} a^{\vee}_{i}\alpha^{\vee}_{i} \in \Fh
\end{equation}
be the null root and the canonical central element of $\Fg$. 
We denote by $(\cdot\,,\,\cdot)$ 
the bilinear form on $\Fh^{\ast}$, which is normalized 
by: $a_{i}^{\vee}=\frac{(\alpha_{i},\alpha_{i})}{2}a_{i}$ 
for all $i \in I$. 
Set $\Fh^{\ast}_{0}:=\bigoplus_{i \in I} 
\BQ \alpha_{i} \subset \Fh^{\ast}$, and let 
$\cl:\Fh^{\ast}_{0} \twoheadrightarrow 
\Fh^{\ast}_{0}/\BQ\delta$ the canonical map
from $\Fh^{\ast}_{0}$ onto the quotient space 
$\Fh^{\ast}_{0}/\BQ\delta$. 
We have a bilinear form (denoted also by $(\cdot\,,\,\cdot)$) 
on $\Fh^{\ast}_{0}/\BQ\delta$ induced 
from the bilinear form $(\cdot\,,\,\cdot)$, 
which is positive-definite. 

We take (and fix) a special vertex 
$0 \in I$ as in \cite[\S5.2]{Kas5}, 
and set $I_{0}:=I \setminus \{0\}$. 
For $i \in I_{0}$, let $\vpi_{i}$ be 
a unique element in $\bigoplus_{i \in I_{0}}
\BQ \alpha_{i}$ such that $\vpi_{i}(\alpha^{\vee}_{j})=
\delta_{i,j}$ for all $j \in I_{0}$. Notice that 
$\Lambda_{i}:=\vpi_{i}+a_{i}^{\vee}\Lambda_{0}$ is an 
$i$-th fundamental weight for $\Fg$, where $\Lambda_{0}$ is 
a $0$-th fundamental weight for $\Fg$. So we may assume that 
all $\vpi_{i}$'s are contained in the integral weight 
lattice $P$. 

%
\subsection{Crystal bases.}
\label{subsec:crystal}

Let $\CB(\infty)$ be the crystal base of the negative 
part $U_{q}^{-}(\Fg)$ with $u_{\infty}$ the 
highest weight element. Denote by 
$e_{i}$ and $f_{i}$ the raising and lowering 
Kashiwara operator on $\CB(\infty)$, respectively, 
and define
$\ve_{i}:\CB(\infty) \rightarrow \BZ$ and 
$\vp_{i}:\CB(\infty) \rightarrow \BZ$ by
\begin{equation}
\ve_{i}(b):=\max\bigl\{n \ge 0 \mid 
e_{i}^{n}b \ne 0\bigr\}, \quad 
\vp_{i}(b):=\ve_{i}(b)+(\wt(b))(\alpha^{\vee}_{i}).
\end{equation}
Denote by $\ast:\CB(\infty) \rightarrow \CB(\infty)$ 
the so-called $\ast$-operation on $\CB(\infty)$ (cf. 
\cite[Theorem 2.1.1]{Kas1} and \cite[\S8.3]{Kas3}). 
We put $e_{i}^{\ast}:=\ast \circ e_{i} \circ \ast$ 
and $f_{i}^{\ast}:=\ast \circ f_{i} \circ \ast$ 
for each $i \in I$. 

\begin{thm}[{\rm cf. \cite[Theorem 2.2.1]{Kas1}}\bf]
\label{thm:kas_emb}
For each $i \in I$, there exists an embedding 
$\Psi_{i}^{-}:\CB(\infty) \hookrightarrow 
\CB(\infty) \otimes \CB_{i}$ of crystals that maps 
$u_{\infty}$ to $u_{\infty} \otimes b_{i}(0)$, 
where $\CB_{i}:=\bigl\{b_{i}(n) \mid n \in \BZ\bigr\}$ 
is a crystal in {\rm \cite[Example 1.2.6]{Kas1}}. 
In addition, if $b=(f_{i}^{\ast})^{k}b_{0}$ 
with some $k \in \BZ_{\ge 0}$
and $b_{0} \in \CB(\infty)$ such that 
$e_{i}^{\ast}b_{0}=0$, then $\Psi^{-}_{i}(b)=b_{0} 
\otimes b_{i}(-k)$. 
\end{thm}

We denote by $\CB(-\infty)$ 
the crystal base of the positive part $U_{q}^{+}(\Fg)$
with $u_{-\infty}$ the lowest weight vector, and 
by $e_{i}$ and $f_{i}$ the raising and lowering 
Kashiwara operator on $\CB(-\infty)$, respectively. We set 
\begin{equation}
\ve_{i}(b):=\vp_{i}(b)-(\wt(b))(\alpha^{\vee}_{i}), \quad 
\vp_{i}(b):=\max\bigl\{n \ge 0 \mid 
f_{i}^{n}b \ne 0\bigr\}.
\end{equation}
We also have the $\ast$-operation 
$\ast:\CB(-\infty) \rightarrow \CB(-\infty)$ 
on $\CB(-\infty)$. We can easily show that 
there exists an embedding 
$\Psi^{+}_{i}:\CB(-\infty) \hookrightarrow \CB_{i} 
\otimes \CB(-\infty)$ of crystals with 
properties similar to $\Psi^{-}_{i}$ 
in Theorem~\ref{thm:kas_emb}.

Let $\CB(\ti{U}_{q}(\Fg))=\bigsqcup_{\lambda \in P}
\CB(U_{q}(\Fg)a_{\lambda})$ be the crystal base of 
the modified quantized universal enveloping algebra 
$\ti{U}_{q}(\Fg)$ with $u_{\lambda}$ the element of 
$\CB(U_{q}(\Fg)a_{\lambda})$ corresponding to 
$a_{\lambda} \in U_{q}(\Fg)a_{\lambda}$
(cf. \cite[Theorem 2.1.2]{Kas2}). 
We denote by $e_{i}$ and $f_{i}$ 
the raising and lowering Kashiwara operator on 
$\CB(\ti{U}_{q}(\Fg))$, and define
$\ve_{i}:\CB(\ti{U}_{q}(\Fg))
\rightarrow \BZ$ and 
$\vp_{i}:\CB(\ti{U}_{q}(\Fg))\rightarrow \BZ$ by 
\begin{equation}
\ve_{i}(b):=\max\bigl\{n \ge 0 \mid 
e_{i}^{n}b \ne 0\bigr\}, 
\qquad 
\vp_{i}(b):=\max\bigl\{n \ge 0 \mid 
f_{i}^{n}b \ne 0\bigr\}. 
\end{equation}
We know the following theorem from \cite[Theorem 3.1.1]{Kas2}.

\begin{thm} \label{thm:cry_mod}
There exists an isomorphism 
$\Xi_{\lambda}:
\CB(U_{q}(\Fg)a_{\lambda}) \stackrel{\sim}{\rightarrow} 
\CB(\infty) \otimes \CT_{\lambda} \otimes \CB(-\infty)$
of crystals such that 
$\Xi_{\lambda}(u_{\lambda})=u_{\infty} \otimes 
t_{\lambda} \otimes u_{-\infty}$, where 
$\CT_{\lambda}:=\bigl\{t_{\lambda}\bigr\}$ 
is a crystal consisting of a 
single element $t_{\lambda}$ of weight $\lambda$ 
{\rm(cf. \cite[Example 7.3]{Kas3})}. 
\end{thm}

We also denote by $\ast:\CB(\ti{U}_{q}(\Fg)) \rightarrow 
\CB(\ti{U}_{q}(\Fg))$ the $\ast$-operation on 
$\CB(\ti{U}_{q}(\Fg))$ (cf. \cite[Theorem 4.3.2]{Kas2}). 
We know the following theorem from \cite[Corollary 4.3.3]{Kas2}. 

\begin{thm} \label{thm:ast}
Let $b \in \CB(U_{q}(\Fg)a_{\lambda})$, and assume that 
$\Xi_{\lambda}(b)=b_{1} \otimes t_{\lambda} \otimes b_{2}$ 
with $b_{1} \in \CB(\infty)$ and $b_{2} \in \CB(-\infty)$. 
Then, $b^{\ast}$ is contained in 
$\CB(U_{q}(\Fg)a_{\lambda^{\prime}})$, where 
$\lambda^{\prime}:=-\lambda-\wt(b_{1})-\wt(b_{2})$, and 
$\Xi_{\lambda^{\prime}}(b^{\ast})=b_{1}^{\ast} \otimes 
t_{\lambda^{\prime}} \otimes b_{2}^{\ast}$. 
\end{thm}

%
\subsection{Crystal base of an extremal weight module.}
\label{subsec:cry_ext}

Since $\CB(\ti{U}_{q}(\Fg))$ is a normal crystal, we 
can define an action of the Weyl group $W$ 
on $\CB(\ti{U}_{q}(\Fg))$ (see \cite[\S7.1]{Kas2});  
for $i \in I$, we define an action of 
simple reflection $r_{i}$ by
\begin{equation} \label{eq:weyl}
r_{i}b:=\begin{cases}
f_{i}^{n}b & \text{if \ } 
n:=(\wt(b))(\alpha^{\vee}_{i}) \ge 0 \\[1.5mm]
e_{i}^{-n}b & \text{if \ } 
n:=(\wt(b))(\alpha^{\vee}_{i}) \le 0. 
\end{cases} \qquad \text{for \ } 
b \in \CB(\ti{U}_{q}(\Fg)). 
\end{equation}
An element $b \in \CB(\ti{U}_{q}(\Fg))$ is called extremal if 
the elements $\bigl\{wb\bigr\}_{w \in W} 
\subset \CB(\ti{U}_{q}(\Fg))$ 
satisfy the following 
condition for all $i \in I$: 
\begin{equation}
\begin{array}{l}
\text{if $\bigl(\wt(wb)\bigr)(\alpha^{\vee}_{i}) \ge 0$, 
then $e_{i}(wb)=0$,} \\[1.5mm]
\text{and if $\bigl(\wt(wb)\bigr)(\alpha^{\vee}_{i}) \le 0$, 
then $f_{i}(wb)=0$.} 
\end{array}
\end{equation}
For $\lambda \in P$, we define a subcrystal 
$\CB(\lambda)$ of $\CB(U_{q}(\Fg)a_{\lambda})$ by
\begin{equation}
\CB(\lambda):=\bigl\{
b \in \CB(U_{q}(\Fg)a_{\lambda}) \mid b^{\ast} \text{\, is 
extremal}\bigr\}. 
\end{equation}
Remark that $u_{\lambda} \in \CB(U_{q}(\Fg)a_{\lambda})$ is 
contained in $\CB(\lambda)$. 
We know from \cite[Proposition 8.2.2]{Kas2} and 
\cite[\S3.1]{Kas5} that $\CB(\lambda)$ is the crystal base
of the extremal weight module $V(\lambda)$ 
of extremal weight $\lambda$ over $U_{q}(\Fg)$. 

%
\section{Some Tools for Crystal Bases.}
\label{sec:pre01}


\subsection{Multiple maps.}
\label{subsec:mult}

We know the following theorem 
(see \cite[Theorem 3.2]{Kas4}). 

\begin{thm} \label{thm:m_inf}
Let $m \in \BZ_{> 0}$. There exists a unique injective map
$S_{m,\infty}:\CB(\infty) \hookrightarrow \CB(\infty)$
such that for all $b \in \CB(\infty)$ and $i \in I$, 
we have
\begin{align}
& \wt(S_{m,\infty}(b))=m \wt(b), \quad
\ve_{i}(S_{m,\infty}(b))=m \ve_{i}(b), \quad 
\vp_{i}(S_{m,\infty}(b))=m \vp_{i}(b), \label{eq:mult01} \\
&
S_{m,\infty}(u_{\infty})=u_{\infty}, \quad
S_{m,\infty}(e_{i}b)=e_{i}^{m}S_{m,\infty}(b), \quad 
S_{m,\infty}(f_{i}b)=f_{i}^{m}S_{m,\infty}(b). \label{eq:mult02}
\end{align}
\end{thm}

\begin{prop} \label{prop:ast_inf}
We set $S_{m,\infty}^{\ast}:=
\ast \circ S_{m,\infty} \circ \ast$. 
Then we have $S_{m,\infty}^{\ast}=S_{m,\infty}$ 
on $\CB(\infty)$. Namely, the $\ast$-operation 
commutes with the map $S_{m,\infty}:\CB(\infty)
\hookrightarrow \CB(\infty)$. 
\end{prop}

The proposition above can be shown in a way similar to
\cite[Theorem 2.3.1]{NS2}. Before giving a proof of the 
proposition, we show the following lemma. 

\begin{lem} \label{lem:com}
The following diagram is commutative\,{\rm:}
\begin{equation}
\begin{CD}
\CB(\infty) @>{\Psi^{-}_{j}}>> 
\CB(\infty) \otimes \CB_{j}\phantom{.} \\
@V{S_{m,\infty}^{\ast}}VV 
@VV{S_{m,\infty}^{\ast} \otimes \CS_{m,j}}V \\
\CB(\infty) @>{\Psi^{-}_{j}}>> 
\CB(\infty) \otimes \CB_{j}.
\end{CD}
\end{equation}
Here $\CS_{m,j}:\CB_{j} \rightarrow \CB_{j}$ is a map defined by 
$\CS_{m,j}(b_{j}(n)):=b_{j}(mn)$. 
\end{lem}

\begin{proof}
For $b \in \CB(\infty)$, there exists 
$b_{0} \in \CB(\infty)$ such that 
$b=(f_{j}^{\ast})^{k}b_{0}$ for some $k \in \BZ_{\ge 0}$
and $e_{j}^{\ast}b_{0}=0$. Then, by Theorem~\ref{thm:kas_emb}, 
we have $\Psi^{-}_{j}(b)=b_{0} \otimes b_{j}(-k)$, and hence
\begin{equation*}
(S^{\ast}_{\infty} \otimes \CS_{m,j})(\Psi^{-}_{j}(b))=
S^{\ast}_{\infty}(b_{0}) \otimes b_{j}(-mk).
\end{equation*}

On the other hand, we see that 
$S_{m,\infty}^{\ast}(b)=(f_{j}^{\ast})^{mk}
S_{m,\infty}^{\ast}(b_{0})$. If 
$e^{\ast}_{j}S_{m,\infty}^{\ast}(b_{0}) \ne 0$, then 
we have $\ve_{j}(S_{m,\infty}(b^{\ast}_{0})) \ge 1$. Since 
$\ve_{j}(S_{m,\infty}(b))=m\ve_{j}(b) \in m\BZ$ for 
all $b \in \CB(\infty)$, we deduce that 
$\ve_{j}(S_{m,\infty}(b^{\ast}_{0})) \ge m$, and hence 
$(e^{\ast}_{j})^{m}S_{m,\infty}^{\ast}(b_{0}) \ne 0$. However, 
since $e^{\ast}_{j}b_{0}=0$, we get
$(e^{\ast}_{j})^{m}S_{m,\infty}^{\ast}(b_{0})
=S_{m,\infty}^{\ast}(e^{\ast}_{j}b_{0})=0$, which is a 
contradiction. Therefore, we conclude that 
$e^{\ast}_{j}S_{m,\infty}^{\ast}(b_{0}) = 0$. 
It follows from Theorem~\ref{thm:kas_emb} that
\begin{equation*}
\Psi^{-}_{j}(S_{m,\infty}^{\ast}(b))=\Psi^{-}_{j}((f_{j}^{\ast})^{mk}
S_{m,\infty}^{\ast}(b_{0}))=
S_{m,\infty}^{\ast}(b_{0}) \otimes b_{j}(-mk). 
\end{equation*}
Hence we have $(S_{m,\infty}^{\ast} \otimes \CS_{m,j})(\Psi^{-}_{j}(b))
=\Psi^{-}_{j}(S_{m,\infty}^{\ast}(b))$. This completes 
the proof of the lemma. 
\end{proof}

\begin{proof}[Proof of Proposition~\ref{prop:ast_inf}]
We will prove that $S^{\ast}_{\infty}(b)=S_{m,\infty}(b)$ for 
$b \in \CB(\infty)_{-\xi}$ by induction on the height $\Ht(\xi)$ of 
$\xi$ (note that $-\wt(b) \in \sum_{i \in I} \BZ_{\ge 0}\alpha_{i}$ for 
all $b \in \CB(\infty)$). If $\Ht(\xi)=0$, then $b$ is 
the highest weight element $u_{\infty} \in \CB(\infty)$, 
and hence the assertion is obvious. 

Assume that $\Ht(\xi) \ge 1$. Then there exists 
some $i \in I$ such that $b_{1}:=e_{i}b \ne 0$. 
If $e^{\ast}_{j}b_{1}=0$ for all 
$j \in I$, then $b_{1}=u_{\infty}$, 
and hence $b=f_{i}u_{\infty}$. Because 
$f^{k}_{i}u_{\infty}$ is a unique element of 
weight $-k\alpha_{i}$ for all $k \in \BZ_{\ge 0}$, and
$\wt(b^{\ast})=\wt(b)$ for all $b \in \CB(\infty)$, 
we deduce that $b^{\ast}=b$, and hence that
\begin{equation*}
S_{m,\infty}^{\ast}(b)=\bigl(S_{m,\infty}(b^{\ast})\bigr)^{\ast}=
\bigl(S_{m,\infty}(b)\bigr)^{\ast}=
\bigl(f_{i}^{m}u_{\infty}\bigr)^{\ast}
=f_{i}^{m}u_{\infty}=S_{m,\infty}(b).
\end{equation*}
So we may assume that there exists $j \in I$
such that $e^{\ast}_{j}b_{1} \ne 0$. Let 
$b_{2} \in \CB(\infty)$ be such that 
$e^{\ast}_{j}b_{2} = 0$ and $b_{1}=
(f_{j}^{\ast})^{k}b_{2}$ for some $k \in \BZ_{\ge 1}$.
Namely, $b=f_{i}(f_{j}^{\ast})^{k}b_{2}$ for some 
$k \ge 1$ and $b_{2} \in \CB(\infty)$ such that 
$e^{\ast}_{j}b_{2} = 0$. 

\vspace{1.5mm}

\noindent {\it Case 1} : $i \ne j$. \, 
We show that $\Psi^{-}_{j}(S_{m,\infty}^{\ast}(b))=
\Psi^{-}_{j}(S_{m,\infty}(b))$ (recall that 
$\Psi^{-}_{j}:\CB(\infty) \hookrightarrow \CB(\infty)
\otimes \CB_{j}$ is an embedding of crystals). 
We have 
\begin{align*}
\Psi^{-}_{j}(b) & =\Psi^{-}_{j}(f_{i}(f_{j}^{\ast})^{k}b_{2})
=f_{i}\Psi^{-}_{j}((f_{j}^{\ast})^{k}b_{2})
=f_{i}(b_{2} \otimes b_{j}(-k)) \\
& = f_{i}b_{2} \otimes b_{j}(-k).
\end{align*}
Here the last equality immediately follows from the 
definition of the tensor product of crystals (see, 
for example, \cite[\S7.3]{Kas3}) and 
the condition that $i \ne j$. Therefore, we obtain
\begin{align*}
\Psi^{-}_{j}(S_{m,\infty}^{\ast}(b)) &=
(S_{m,\infty}^{\ast} \otimes \CS_{m,j})(\Psi^{-}_{j}(b)) 
\quad \text{by Lemma~\ref{lem:com}} \\
& =S_{m,\infty}^{\ast}(f_{i}b_{2}) \otimes b_{j}(-mk) \\
& =S_{m,\infty}(f_{i}b_{2}) \otimes b_{j}(-mk)
\quad \text{by the inductive assumption} \\
& =f_{i}^{m}S_{m,\infty}(b_{2}) \otimes b_{j}(-mk).
\end{align*}
On the other hand, 
\begin{align*}
S_{m,\infty}(b) & 
  =S_{m,\infty}(f_{i}(f_{j}^{\ast})^{k}b_{2})
  =f_{i}^{m}S_{m,\infty}((f_{j}^{\ast})^{k}b_{2}) \\
& =f_{i}^{m}(f_{j}^{\ast})^{mk}S_{m,\infty}(b_{2})
\quad \text{by the inductive assumption.}
\end{align*}
As in the proof of Lemma~\ref{lem:com}, 
we deduce that 
$e_{j}^{\ast}S_{m,\infty}^{\ast}(b_{2})=0$, 
and hence 
$e_{j}^{\ast}S_{m,\infty}(b_{2})=
e_{j}^{\ast}S_{m,\infty}^{\ast}(b_{2})=0$ 
by the inductive assumption. Therefore, 
\begin{align*}
\Psi^{-}_{j}(S_{m,\infty}(b))  
& = \Psi^{-}_{j}(f_{i}^{m}(f_{j}^{\ast})^{mk}
S_{m,\infty}(b_{2}))
= f_{i}^{m}\Psi^{-}_{j}((f_{j}^{\ast})^{mk}
S_{m,\infty}(b_{2})) \\
& =f_{i}^{m}(S_{m,\infty}(b_{2}) \otimes b_{j}(-mk))
=(f_{i}^{m}S_{m,\infty}(b_{2})) \otimes b_{j}(-mk).
\end{align*}
Here the last equality immediately follows again 
from the definition of the tensor product of crystals 
and the condition that $i \ne j$. 
Thus, we get that $\Psi^{-}_{j}(S_{m,\infty}^{\ast}(b))=
\Psi^{-}_{j}(S_{m,\infty}(b))$, and hence 
$S_{m,\infty}^{\ast}(b)=S_{m,\infty}(b)$. 

\vspace{1.5mm}

\noindent {\it Case 2} : $i=j$. \, 
As in Case 1, we have 
$\Psi^{-}_{j}(b)=f_{i}(b_{2} \otimes b_{i}(-k))$. 
We deduce from the definition of the tensor product of 
crystals that 
\begin{equation*}
\Psi^{-}_{i}(b)=f_{i}(b_{2} \otimes b_{i}(-k))
=\begin{cases}
f_{i}b_{2} \otimes b_{i}(-k) 
& \text{if } \vp_{i}(b_{2}) > k, \\[1.5mm]
b_{2} \otimes b_{i}(-k-1) 
& \text{if } \vp_{i}(b_{2}) \le k.
\end{cases}
\end{equation*}
Hence, as in Case 1, we get 
\begin{equation*}
\Psi^{-}_{i}(S_{m,\infty}^{\ast}(b))=
\begin{cases}
f_{i}^{m}S_{m,\infty}(b_{2}) \otimes b_{i}(-mk) 
& \text{if } \vp_{i}(b_{2}) > k, \\[1.5mm]
S_{m,\infty}(b_{2}) \otimes b_{i}(-mk-m) 
& \text{if } \vp_{i}(b_{2}) \le k.
\end{cases}
\end{equation*}
On the other hand, in exactly the same way as in Case 1, 
we can show that $\Psi^{-}_{i}(S_{m,\infty}(b))  
 =f_{i}^{m}(S_{m,\infty}(b_{2}) \otimes b_{i}(-mk))$. 
Because $\vp_{i}(S_{m,\infty}(b_{2}))=m\vp_{i}(b_{2})$ 
by \eqref{eq:mult01}, we deduce 
from the definition of the tensor product of 
crystals that 
\begin{equation*} 
f_{i}^{m}(S_{m,\infty}(b_{2}) \otimes b_{i}(-mk))
=\begin{cases}
f_{i}^{m}S_{m,\infty}(b_{2}) \otimes b_{i}(-mk) 
& \text{if } \vp_{i}(b_{2}) > k, \\[1.5mm]
S_{m,\infty}(b_{2}) \otimes b_{i}(-mk-m) 
& \text{if } \vp_{i}(b_{2}) \le k.
\end{cases}
\end{equation*}
Therefore we obtain that 
$\Psi^{-}_{i}(S_{m,\infty}^{\ast}(b))=\Psi^{-}_{i}(S_{m,\infty}(b))$, 
and hence $S_{m,\infty}^{\ast}(b)=S_{m,\infty}(b)$. 
Thus, we have proved the proposition. 
\end{proof}

\begin{rem} \label{rem:ast_inf}
A similar result holds for the crystal base 
$\CB(-\infty)$. Namely, for each $m \in \BZ_{>0}$, 
there exists a unique injective map $S_{m,-\infty}:
\CB(-\infty) \hookrightarrow \CB(-\infty)$ with properties 
similar to $S_{m,\infty}$ in Theorem~\ref{thm:m_inf}, 
and it commutes with the $\ast$-operation on $\CB(-\infty)$. 
\end{rem}

For $m \in \BZ_{>0}$, we define an injective map 
$\ti{S}_{m,\lambda}:\CB(U_{q}(\Fg)a_{\lambda}) 
\hookrightarrow \CB(U_{q}(\Fg)a_{m\lambda})$ 
as in the following commutative diagram 
(cf. Theorem~\ref{thm:cry_mod}):
\begin{equation}
\begin{CD}
\CB(U_{q}(\Fg)a_{\lambda}) @>{\Xi_{\lambda}}>{\sim}> 
\CB(\infty) \otimes \CT_{\lambda} \otimes 
\CB(-\infty)\phantom{,} \\
@V{\ti{S}_{m,\lambda}}VV 
@VV{S_{m,\infty} \otimes \tau_{m,\lambda} \otimes S_{m,-\infty}}V \\
\CB(U_{q}(\Fg)a_{m\lambda}) @<{\Xi^{-1}_{m\lambda}}<{\sim}< 
\CB(\infty) \otimes \CT_{m\lambda} \otimes \CB(-\infty),
\end{CD}
\end{equation}
where $\tau_{m,\lambda}:\CT_{\lambda} \rightarrow 
\CT_{m\lambda}$ is defined by $\tau_{m,\lambda}
(t_{\lambda}):=t_{m\lambda}$.
We define $\ti{S}_{m}:\ti{U}_{q}(\Fg) 
\hookrightarrow \ti{U}_{q}(\Fg)$
as the direct sum of all $\ti{S}_{m,\lambda}$'s.

\begin{prop} \label{prop:m_mod}
The maps $\ti{S}_{m,\lambda}:\CB(U_{q}(\Fg)a_{\lambda}) 
\hookrightarrow \CB(U_{q}(\Fg)a_{m\lambda})$ and 
$\ti{S}_{m}:\CB(\ti{U}_{q}(\Fg)) 
\hookrightarrow \CB(\ti{U}_{q}(\Fg))$ have 
properties similar to $S_{m,\infty}$ in Theorem~\ref{thm:m_inf}. 
In addition, the map $\ti{S}_{m}$ commutes with the 
$\ast$-operation on $\CB(\ti{U}_{q}(\Fg))$. 
\end{prop}

\begin{proof}
The first assertion immediately follows from 
Theorem~\ref{thm:m_inf}, Remark~\ref{rem:ast_inf}, 
and the definition of the tensor product of crystals 
(see also \cite[Appendix B]{Kas5}). 
Let us show the second assertion. We set $\ti{S}_{m}^{\ast}:=
\ast \circ \ti{S}_{m} \circ \ast$. It suffices to 
prove the following claim: 

\vspace{1.5mm}

\noindent
{\it Claim.\ } Let $\lambda \in P$, and 
$b \in \CB(U_{q}(\Fg)a_{\lambda})$. Then we have 
that $\ti{S}_{m}^{\ast}(b) \in \CB(U_{q}(\Fg)a_{m\lambda})$, 
and that $\Xi_{m\lambda}(\ti{S}_{m}^{\ast}(b))=
\Xi_{m\lambda}(\ti{S}_{m}(b))$. 

\vspace{1.5mm}

\noindent Assume that 
$\Xi_{\lambda}(b)=b_{1} \otimes t_{\lambda} \otimes 
b_{2}$ with $b_{1} \in \CB(\infty)$ and 
$b_{2} \in \CB(-\infty)$. Then we see by the definition of 
$\ti{S}_{m}$ that 
\begin{equation*}
\Xi_{m\lambda}(\ti{S}_{m}(b))=(S_{m,\infty} \otimes \tau_{m,\lambda}
\otimes S_{m,-\infty})(\Xi_{\lambda}(b))=
S_{m,\infty}(b_{1}) \otimes t_{m\lambda} \otimes 
S_{m,-\infty}(b_{2}). 
\end{equation*}
On the other hand, we know from Theorem~\ref{thm:ast}
that $b^{\ast} \in \CB(U_{q}(\Fg)a_{\lambda^{\prime}})$ and 
$\Xi_{\lambda^{\prime}}(b^{\ast})=b_{1}^{\ast} \otimes 
t_{\lambda^{\prime}} \otimes b_{2}^{\ast}$, where 
$\lambda^{\prime}:=-\lambda-\wt(b_{1})-\wt(b_{2})$. 
Hence we have 
\begin{equation*}
\Xi_{m\lambda^{\prime}}(\ti{S}_{m}(b^{\ast}))=
(S_{m,\infty} \otimes \tau_{m,\lambda}
\otimes S_{m,-\infty})(\Xi_{\lambda^{\prime}}(b^{\ast}))=
S_{m,\infty}(b_{1}^{\ast}) \otimes t_{m\lambda^{\prime}} \otimes 
S_{m,-\infty}(b_{2}^{\ast}). 
\end{equation*}
We deduce again from Theorem~\ref{thm:ast}
that $\ti{S}_{m}^{\ast}(b)=
\bigl(\ti{S}_{m}(b^{\ast})\bigr)^{\ast} \in 
\CB(U_{q}(\Fg)a_{m\lambda})$, and that 
\begin{align*}
& \Xi_{m\lambda}(\ti{S}_{m}^{\ast}(b)) 
=S_{m,\infty}^{\ast}(b_{1}) 
\otimes t_{m\lambda} \otimes 
S_{m,-\infty}^{\ast}(b_{2}) \\
& \hspace{5mm} 
=S_{m,\infty}(b_{1}) 
\otimes t_{m\lambda} \otimes 
S_{m,-\infty}(b_{2}) \quad 
\text{by Proposition~\ref{prop:ast_inf} and 
Remark~\ref{rem:ast_inf}}.
\end{align*}
Thus, we obtain $\Xi_{m\lambda}(\ti{S}_{m}^{\ast}(b))=
\Xi_{m\lambda}(\ti{S}_{m}(b))$, as desired. 
\end{proof}

\begin{thm}\label{thm:m_ext}
Let $m \in \BZ_{> 0}$. There exists an injective map 
$S_{m,\lambda}:\CB(\lambda) \hookrightarrow \CB(m\lambda)$ 
such that $S_{m,\lambda}(u_{\lambda})=u_{m\lambda}$ 
and such that for all $b \in \CB(\infty)$ and 
$i \in I$, we have
\begin{align}
& \wt(S_{m,\lambda}(b))=m \wt(b), \quad
\ve_{i}(S_{m,\lambda}(b))=m \ve_{i}(b), \quad 
\vp_{i}(S_{m,\lambda}(b))=m \vp_{i}(b), \\
&
S_{m,\lambda}(e_{i}b)=e_{i}^{m}S_{m,\lambda}(b), \quad 
S_{m,\lambda}(f_{i}b)=f_{i}^{m}S_{m,\lambda}(b).
\end{align}
\end{thm}

\begin{proof}
Set $S_{m,\lambda}:=\ti{S}_{m}|_{\CB(\lambda)}$. Then 
it is obvious from Proposition~\ref{prop:m_mod} that 
$S_{m,\lambda}(\CB(\lambda)) \subset \CB(U_{q}(\Fg)a_{m\lambda})$. 
So we need only show that $\bigl(S_{m,\lambda}(b)\bigr)^{\ast}$
is extremal for every $b \in \CB(\lambda)$. We can easily check 
that the action of the Weyl group $W$ commutes with $S_{m,\lambda}$. So
it follows from Proposition~\ref{prop:m_mod} that
\begin{equation*}
w\bigl((S_{m,\lambda}(b))^{\ast}\bigr)=w S_{m,\lambda}(b^{\ast})
=S_{m,\lambda}(wb^{\ast}) \quad 
\text{for all \ } b \in \CB(\lambda)
\text{\ and\ } w \in W. 
\end{equation*}
Assume that $\wt(b^{\ast})=\mu$. Then we see that
$\wt\bigl((S_{m,\lambda}(b))^{\ast}\bigr)=m\mu$.
Suppose that $(w(m\mu))(\alpha^{\vee}_{i}) \ge 0$ and 
$e_{i}\bigl(w((S_{m,\lambda}(b))^{\ast})\bigr) \ne 0$. 
As in the proof of Lemma~\ref{lem:com}, 
we deduce that $e_{i}^{m}
\bigl(w((S_{m,\lambda}(b))^{\ast})\bigr) \ne 0$. 
Hence we have
\begin{equation*}
S_{m,\lambda}\bigl(e_{i}(wb^{\ast})\bigr)=
e_{i}^{m}S_{m,\lambda}(wb^{\ast})=
e_{i}^{m}\bigl(w S_{m,\lambda}(b^{\ast})\bigr)=
e_{i}^{m}\bigl(w((S_{m,\lambda}(b))^{\ast})\bigr) \ne 0.
\end{equation*}
However, 
since $(w(\mu))(\alpha^{\vee}_{i}) \ge 0$ and 
$b^{\ast}$ is extremal, 
we have $e_{i}(wb^{\ast})=0$, and hence 
$S_{m,\lambda}\bigl(e_{i}(wb^{\ast})\bigr)=0$, 
which is a contradiction.
Therefore, we get that 
$e_{i}\bigl(w((S_{m,\lambda}(b))^{\ast})\bigr)=0$. 
Similarly, we can prove that if $(w(m\mu))
(\alpha^{\vee}_{i}) \le 0$, then 
$f_{i}\bigl(w((S_{m,\lambda}(b))^{\ast})\bigr)=0$. 
This completes the proof of the theorem.
\end{proof}


\subsection{Embedding into tensor products.}
\label{subsec:fact01}

In this subsection, 
we assume that $\Fg$ is an affine Lie algebra
(for the notation, see \S\ref{subsec:affine}).
We know the following theorem from \cite{B}, \cite{N}
in the symmetric case, and from \cite{BN} in the 
nonsymmetric case. 
\begin{thm} \label{thm:emb}
Let $\lambda:=\sum_{i \in I_{0}}m_{i}\vpi_{i}$ with 
$m_{i} \in \BZ_{\ge 0}$, and let $\CB_{0}(\lambda)$ 
be the connected component of $\CB(\lambda)$ containing 
the extremal weight element $u_{\lambda}$ of weight $\lambda$.
There exists a unique embedding 
$\CB_{0}(\lambda) \hookrightarrow \bigotimes_{i \in I_{0}}
\CB(\vpi_{i})^{\otimes m_{i}}$ of crystals that maps 
$u_{\lambda}$ to $\bigotimes_{i \in I_{0}}
u_{\vpi_{i}}^{\otimes m_{i}}$. In particular, 
we have an embedding 
\begin{equation}
G_{m,\vpi_{i}}:\CB_{0}(m\vpi_{i}) \hookrightarrow 
\CB(\vpi_{i})^{\otimes m}
\end{equation}
of crystals that maps $u_{m\vpi_{i}}$ to 
$u_{\vpi_{i}}^{\otimes m}$. 
\end{thm}

Since $\CB(\vpi_{i})$ is connected (see \cite[Theorem 5.5]{Kas5}), 
we see that $S_{m,\vpi_{i}}(\CB(\vpi_{i})) \subset 
\CB_{0}(m\vpi_{i})$. Hence we can define 
$\sigma_{m,\vpi_{i}}:\CB(\vpi_{i}) 
\hookrightarrow \CB(\vpi_{i})^{\otimes m}$ by 
$\sigma_{m,\vpi_{i}}:=G_{m,\vpi_{i}} \circ S_{m,\vpi_{i}}$ 
for each $m \in \BZ_{> 0}$. Remark that $\sigma_{m,\vpi_{i}}$ 
has the following properties: 
\begin{align}
& 
\wt(\sigma_{m,\vpi_{i}}(b))=m \wt(b), \quad
\ve_{j}(\sigma_{m,\vpi_{i}}(b))=m \ve_{j}(b), \quad 
\vp_{j}(\sigma_{m,\vpi_{i}}(b))=m \vp_{j}(b), 
\label{eq:mult03} \\
&
\sigma_{m,\vpi_{i}}(u_{\vpi_{i}})=
u_{\vpi_{i}}^{\otimes m}, \quad
\sigma_{m,\vpi_{i}}(e_{j}b)=
e_{j}^{m}\sigma_{m,\vpi_{i}}(b), \quad 
\sigma_{m,\vpi_{i}}(f_{j}b)=
f_{j}^{m}\sigma_{m,\vpi_{i}}(b). \label{eq:mult04}
\end{align}

\begin{lem} \label{lem:mult_com}
Let $m,\,n \in \BZ_{> 0}$. Then we have 
$\sigma_{mn,\vpi_{i}}=\sigma_{n,\vpi_{i}}^{\otimes m} \circ 
\sigma_{m,\vpi_{i}}$.
\end{lem}

\begin{proof}
Since $\CB(\vpi_{i})$ is connected, 
every $b \in \CB(\vpi_{i})$ is of the form 
$b=x_{j_{1}}x_{j_{2}} \cdots x_{j_{k}}u_{\vpi_{i}}$ for 
some $j_{1},\,j_{2},\,\dots,\,j_{k} \in I$, where 
$x_{j}$ is either $e_{j}$ or $f_{j}$. We will show by 
induction on $k$ that 
$\sigma_{mn,\vpi_{i}}(b)=\sigma_{n,\vpi_{i}}^{\otimes m} \circ 
\sigma_{m,\vpi_{i}}(b)$ for all $b \in \CB(\vpi_{i})$. 
If $k=0$, then the assertion is obvious, 
since $b=u_{\vpi_{i}}$. Assume that $k \ge 1$. We set 
$b^{\prime}:=x_{j_{2}} \cdots x_{j_{k}}u_{\vpi_{i}}$, and 
$\sigma_{m,\vpi_{i}}(b^{\prime})=:u_{1} \otimes u_{2} \otimes 
\cdots \otimes u_{m} \in \CB(\vpi_{i})^{\otimes m}$. Assume that 
\begin{equation*}
\sigma_{m,\vpi_{i}}(b)=
x_{j_{1}}^{m}\sigma_{m,\vpi_{i}}(b^{\prime})=
x_{j_{1}}^{k_{1}}u_{1} \otimes x_{j_{1}}^{k_{2}}u_{2} 
\otimes \cdots \otimes x_{j_{1}}^{k_{m}}u_{m}.
\end{equation*}
for some $k_{1},\,k_{2},\,\dots,\,k_{m} \in \BZ_{\ge 0}$. 
Then we have 
\begin{equation*}
\sigma_{n,\vpi_{i}}^{\otimes m} \circ 
\sigma_{m,\vpi_{i}}(b)=
x_{j_{1}}^{nk_{1}}\sigma_{n,\vpi_{i}}(u_{1}) 
\otimes x_{j_{1}}^{nk_{2}}\sigma_{n,\vpi_{i}}(u_{2})
\otimes \cdots \otimes x_{j_{1}}^{nk_{m}}
\sigma_{n,\vpi_{i}}(u_{m}).
\end{equation*}
Remark (cf. \cite[Lemma 1.3.6]{Kas1}) 
that for all $u_{1} \otimes u_{2} 
\otimes \dots \otimes u_{m} \in 
\CB(\vpi_{i})^{\otimes m}$, 
\begin{equation*}
x_{j}(u_{1} \otimes u_{2} \otimes \cdots \otimes u_{m})=
u_{1} \otimes u_{2} \otimes \cdots \otimes x_{j}u_{l} \otimes 
\cdots \otimes u_{m}
\end{equation*}
if and only if
\begin{align*}
& 
x_{j}^{n}(\sigma_{n,\vpi_{i}}(u_{1}) \otimes 
\sigma_{n,\vpi_{i}}(u_{2}) \otimes \cdots \otimes 
\sigma_{n,\vpi_{i}}(u_{m}))= \\
& \hspace{15mm}
\sigma_{n,\vpi_{i}}(u_{1}) \otimes 
\sigma_{n,\vpi_{i}}(u_{2}) \otimes \cdots \otimes 
x_{j}^{n}\sigma_{n,\vpi_{i}}(u_{l}) \otimes 
\cdots \otimes \sigma_{n,\vpi_{i}}(u_{m}). 
\end{align*}
So we obtain 
\begin{align*}
\sigma_{n,\vpi_{i}}^{\otimes m} \circ 
\sigma_{m,\vpi_{i}}(b)
 & =x_{j_{1}}^{mn}\bigl(
\sigma_{n,\vpi_{i}}(u_{1}) 
\otimes \sigma_{n,\vpi_{i}}(u_{2})
\otimes \cdots \otimes \sigma_{n,\vpi_{i}}(u_{m})\bigr) \\
& =x_{j_{1}}^{mn}\bigl(\sigma_{n,\vpi_{i}}^{\otimes m} \circ 
\sigma_{m,\vpi_{i}}(b^{\prime})\bigr).
\end{align*}
Since $\sigma_{n,\vpi_{i}}^{\otimes m} \circ 
\sigma_{m,\vpi_{i}}(b^{\prime})=
\sigma_{mn,\vpi_{i}}(b^{\prime})$ by the inductive assumption, 
and since $\sigma_{mn,\vpi_{i}}(b)=x_{j_{1}}^{mn}
\sigma_{mn,\vpi_{i}}(b^{\prime})$, we obtain 
$\sigma_{n,\vpi_{i}}^{\otimes m} \circ 
\sigma_{m,\vpi_{i}}(b)=\sigma_{mn,\vpi_{i}}(b)$. 
\end{proof}

For each $w \in W$, we set $u_{w\vpi_{i}}:=wu_{\vpi_{i}} 
\in \CB(\vpi_{i})$. By \cite[Proposition 5.8]{Kas5}, 
we see that $u_{w\lambda}$ is well-defined. 
We can easily show the following lemma. 

\begin{lem} \label{lem:wlambda}
For each $m \in \BZ_{ > 0}$ and $w \in W$, 
we have $\sigma_{m,\vpi_{i}}(u_{w\vpi_{i}})=
(u_{w\vpi_{i}})^{\otimes m}$. 
\end{lem}

\begin{prop} \label{prop:extb}
Let $b \in \CB(\vpi_{i})$. Assume that 
$b=x_{j_{1}}x_{j_{2}} \cdots x_{j_{k}}u_{\vpi_{i}}$, 
where $x_{j}$ is either $e_{j}$ or $f_{j}$,
and set $b_{l}:=x_{j_{l}}x_{j_{l+1}} \cdots x_{i_{k}}
u_{\vpi_{i}}$ for $l=1,\,2\,\dots,\,k+1$ {\rm(}here 
$b_{k+1}:=u_{\vpi_{i}}${\rm)}. Then 
there exists sufficiently large $m \in \BZ$ such that for 
all $l=1,\,2\,\dots,\,k+1$, 
\begin{equation}
\sigma_{m,\vpi_{i}}(b_{l})=u_{w_{l,1}\vpi_{i}} \otimes 
u_{w_{l,2}\vpi_{i}} \otimes \cdots \otimes 
u_{w_{l,m}\vpi_{i}}
\end{equation}
for some $w_{l,1},\,w_{l,2},\,\dots,\,w_{l,m} \in W$. 
\end{prop}

\begin{proof}
We show the assertion by induction 
on $k$. If $k=0$, then the assertion is obvious. 
Assume that $k \ge 1$. By the inductive assumption, 
there exists $m \in \BZ_{>0}$ such that
$\sigma_{m,\vpi_{i}}(b_{l})$ is of the desired form 
for all $l=2,\,\dots,k+1$. Assume that 
\begin{align*}
\sigma_{m,\vpi_{i}}(b_{1})
& =\sigma_{m,\vpi_{i}}(x_{j_{1}}b_{2})
  =x_{j_{1}}^{m}\sigma_{m,\vpi_{i}}(b_{2}) \\
& =x_{j_{1}}^{c_{1}}u_{w_{2,1}\vpi_{i}} \otimes 
x_{j_{1}}^{c_{2}}u_{w_{2,2}\vpi_{i}} \otimes \cdots \otimes 
x_{j_{1}}^{c_{m}}u_{w_{2,m}\vpi_{i}}.
\end{align*}
for some $c_{1},\,c_{2},\,\dots,\,c_{m} \in \BZ_{\ge 0}$.
We can easily check by Lemma~\ref{lem:wlambda} and 
\cite[Lemma 1.3.6]{Kas1} 
that if $n_{p} \in \BZ_{> 0}$ satisfies the condition
$(w_{2,p}\vpi_{i})(\alpha^{\vee}_{j_{1}}) \mid n_{p}c_{p}$, 
then $\sigma_{n_{p},\vpi_{i}}(x_{j_{1}}^{c_{p}}u_{w_{2,p}\vpi_{i}})=
u_{w_{1}\vpi_{i}} \otimes u_{w_{2}\vpi_{i}} 
\otimes \dots \otimes u_{w_{n}\vpi_{i}}$ 
for some $w_{1},\,w_{2},\,\dots,\,w_{n} \in W$. 
Therefore, by Lemma~\ref{lem:wlambda}, 
we see that there exists $N \gg 0$ 
(for example, put $N=\prod_{p=1}^{m}n_{p}$) such that 
\begin{equation*}
(\sigma_{N,\vpi_{i}})^{\otimes m} \circ 
\sigma_{m,\vpi_{i}}(b_{1})=
u_{w_{1,1}\vpi_{i}} \otimes u_{w_{1,2}\vpi_{i}} 
\otimes \dots \otimes u_{w_{1,Nm}\vpi_{i}}
\end{equation*}
for some $w_{1,1},\,w_{1,2},\,\dots,\,
w_{1,Nm} \in W$. Furthermore, we deduce from 
Lemma~\ref{lem:wlambda} that 
$(\sigma_{N,\vpi_{i}})^{\otimes m} \circ 
\sigma_{m,\vpi_{i}}(b_{l})$ is of the desired form 
for all $l=2,\,\dots,\,k+1$. It follows from 
Lemma~\ref{lem:mult_com} that 
$(\sigma_{N,\vpi_{i}})^{\otimes m} \circ 
\sigma_{m,\vpi_{i}}=\sigma_{Nm,\vpi_{i}}$. 
We have thus proved the proposition. 
\end{proof}

%
\section{Preliminary Results.}
\label{sec:pre02}

%
\subsection{Some tools for path models.}
\label{subsec:path}

A path is, by definition, a piecewise linear, continuous 
map $\pi:[0,1] \rightarrow \BQ \otimes_{\BZ} P$ such that $\pi(0)=0$. 
We regard two paths $\pi$ and $\pi^{\prime}$ as equivalent
if there exist piecewise linear, nondecreasing, surjective, 
continuous maps $\psi,\,\psi^{\prime}:[0,1] \rightarrow [0,1]$ 
(reparametrization) such that 
$\pi \circ \psi=\pi^{\prime} \circ \psi$. 
We denote by $\BP$ the set of paths (modulo reparametrization) 
such that $\pi(1) \in P$, 
and by $e_{i}$ and $f_{i}$ the raising and lowering 
root operator (see \cite[\S1]{L2}). By using root operators, 
we can endow $\BP$ with a normal crystal structure 
(see \cite[\S1 and \S2]{L2}); we set $\wt(\pi):=\pi(1)$, and 
define $\ve_{i}:\BP \rightarrow \BZ$ and 
$\vp_{i}:\BP \rightarrow \BZ$ by
\begin{equation}
\ve_{i}(\pi):=\max\bigl\{n \ge 0 \mid 
e_{i}^{n}\pi \ne 0\bigr\}, \qquad 
\vp_{i}(\pi):=\max\bigl\{n \ge 0 \mid 
f_{i}^{n}\pi \ne 0\bigr\}. 
\end{equation}

For an (arbitrary) integral weight $\lambda \in P$, 
we denote by $\BB(\lambda)$ the set of 
Lakshmibai--Seshadri paths of shape $\lambda$ 
(see \cite[\S4]{L2}), and set $\pi_{\lambda}(t):=
t\lambda \in \BB(\lambda)$. Denote by $\BB_{0}(\lambda)$ 
the connected component of $\BB(\lambda)$ containing 
$\pi_{\lambda}$. We obtain the following 
lemma by \cite[Lemma 2.4]{L2}. 

\begin{lem} \label{lem:mp}
For $\pi \in \BP$, we define $S_{m}:\BP \hookrightarrow 
\BP$ by $S_{m}(\pi):=m\pi$, where $(m\pi)(t):=m\pi(t)$ 
for $t \in [0,1]$. Then we have $S_{m}(\BB_{0}(\lambda))=
\BB_{0}(m\lambda)$. In addition, the map $S_{m}$ has 
properties similar to $S_{m,\infty}$ in Theorem~\ref{thm:m_inf}. 
\end{lem}

For paths $\pi_{1},\,\pi_{2} \in \BP$, we define 
a concatenation $\pi_{1} \ast \pi_{2} \in \BP$ 
as in \cite[\S1]{L2}. Because $\pi_{\lambda} \ast 
\pi_{\lambda} \ast \cdots \ast \pi_{\lambda}$
($m$-times) is just $\pi_{m\lambda}$ modulo reparametrization, 
we obtain the following lemma.

\begin{lem} \label{lem:embp}
We have a canonical embedding 
$G_{m,\lambda}:\BB_{0}(m\lambda) \hookrightarrow 
\BB(\lambda)^{\ast m}$ of crystals that maps $\pi_{m\lambda}$
to $\pi_{\lambda}^{\ast m}$, where 
$\BB(\lambda)^{\ast m}:=\bigl\{
\pi_{1} \ast \pi_{2} \ast \cdots \ast \pi_{m} 
\mid \pi_{i} \in \BB(\lambda)\bigr\}$,
and $\pi_{\lambda}^{\ast m}:=\pi_{\lambda} \ast 
\pi_{\lambda} \ast \cdots \ast \pi_{\lambda} \in 
\BB(\lambda)^{\ast m}$.
\end{lem}

By combining Lemmas~\ref{lem:mp} and \ref{lem:embp}, 
we get an embedding $\sigma_{m,\lambda}:\BB_{0}(\lambda) 
\hookrightarrow \BB(\lambda)^{\ast m}$ defined by 
$\sigma_{m,\lambda}:=G_{m,\lambda} \circ S_{m}$. It can easily 
be seen that this map has properties similar to 
\eqref{eq:mult03} and \eqref{eq:mult04}. 

Since $\BB(\lambda)$ is a normal crystal, 
we can define an action of the Weyl group $W$ on 
$\BB(\lambda)$ (cf. \eqref{eq:weyl}; see also 
\cite[Theorem 8.1]{L2}). We set $\pi_{w\lambda}:=
w\pi_{\lambda}$ for $w \in W$. Note that 
$(w\pi_{\lambda})(t)=t(w\lambda)$ for each $w \in W$. 
Using \cite[Lemma 2.7]{L2}, 
we can prove the following proposition 
in a way similar to Proposition~\ref{prop:extb}.

\begin{prop} \label{prop:extp}
Let $\pi \in \BB_{0}(\lambda)$. Assume that 
$\pi=x_{j_{1}}x_{j_{2}} \cdots x_{j_{k}}\pi_{\lambda}$, 
where $x_{j}$ is either $e_{j}$ or $f_{j}$,
and set $\pi_{l}:=x_{j_{l}}x_{j_{l+1}} \cdots x_{i_{k}}
\pi_{\lambda}$ for $l=1,\,2,\,\dots,\,k+1$ {\rm(}here
$\pi_{k+1}:=\pi_{\lambda}${\rm)}. Then 
there exists sufficiently large $m \in \BZ$ such that for 
all $l=1,\,2\,\dots,\,k+1$, 
\begin{equation}
\sigma_{m,\lambda}(\pi_{l})=\pi_{w_{l,1}\lambda} \ast 
\pi_{w_{l,2}\lambda} \ast \cdots \ast 
\pi_{w_{l,m}\lambda}
\end{equation}
for some $w_{l,1},\,w_{l,2},\,\dots,\,w_{l,m} \in W$. 
\end{prop}


\subsection{Preliminary lemmas.}
\label{subsec:fact02}

In this subsection, $\Fg$ is again of affine type 
(for the notation, see \S\ref{subsec:affine}). 
By using \cite[Lemma 2.1 c)]{L2}, we can easily show
the following lemma. 

\begin{lem}
Let $i \in I_{0}$. 
For each $w \in W$ and $j \in I$, we have 
$\wt(\pi_{w\vpi_{i}})=\wt(u_{w\vpi_{i}})$, 
$\ve_{j}(\pi_{w\vpi_{i}})=\ve_{j}(u_{w\vpi_{i}})$, 
and $\vp_{j}(\pi_{w\vpi_{i}})=\vp_{j}(u_{w\vpi_{i}})$.
\end{lem}

It follows from \cite[Lemma 1.3.6]{Kas1}, 
\cite[Lemma 2.7]{L2}, and the lemma above that
\begin{equation*}
x_{j}^{k}(u_{w_{1}\vpi_{i}} \otimes 
u_{w_{2}\vpi_{i}} \otimes \cdots \otimes 
u_{w_{m}\vpi_{i}})=
x_{j}^{k_{1}}u_{w_{1}\vpi_{i}} \otimes 
x_{j}^{k_{2}}u_{w_{2}\vpi_{i}} \otimes \cdots \otimes 
x_{j}^{k_{m}}u_{w_{m}\vpi_{i}}
\end{equation*}
for some $k_{1},\,k_{2},\,\dots,\,k_{m} \in \BZ_{\ge 0}$ 
if and only if 
\begin{equation*}
x_{j}^{k}(\pi_{w_{1}\vpi_{i}} \ast 
\pi_{w_{2}\vpi_{i}} \ast \cdots \ast 
\pi_{w_{m}\vpi_{i}})=
x_{j}^{k_{1}}\pi_{w_{1}\vpi_{i}} \ast 
x_{j}^{k_{2}}\pi_{w_{2}\vpi_{i}} \ast \cdots \ast 
x_{j}^{k_{m}}\pi_{w_{m}\vpi_{i}}
\end{equation*}
for every $k \in \BZ_{\ge 0}$, $m \in \BZ_{> 0}$ and 
$w_{1},\,w_{2},\,\dots,\,w_{m} \in W$. So we obtain 
the following lemma.

\begin{lem} \label{lem:extbp}
{\rm (1)} \ 
Let $b=x_{j_{1}}x_{j_{2}} \cdots x_{j_{k}}u_{\vpi_{i}} 
\in \CB(\vpi_{i})$. Take $m \in \BZ_{> 0}$ such that 
the assertion of Proposition~\ref{prop:extb} holds, 
and assume that 
$\sigma_{m,\vpi_{i}}(b)=u_{w_{1}\vpi_{i}} \otimes 
u_{w_{2}\vpi_{i}} \otimes \cdots \otimes u_{w_{m}\vpi_{i}}$. 
Then we have $\pi:=x_{j_{1}}x_{j_{2}} \cdots x_{j_{k}}
\pi_{\vpi_{i}} \ne 0$, and 
$\sigma_{m,\vpi_{i}}(\pi)=\pi_{w_{1}\vpi_{i}} \ast 
\pi_{w_{2}\vpi_{i}} \ast \cdots \ast \pi_{w_{m}\vpi_{i}}$. 

\vspace{1.5mm}

\noindent {\rm (2)} \ 
The converse of {\rm (1)} holds. Namely, 
let $\pi=x_{j_{1}}x_{j_{2}} \cdots x_{j_{k}}\pi_{\vpi_{i}} 
\in \BB(\vpi_{i})$. Take $m \in \BZ_{> 0}$ such that 
the assertion of Proposition~\ref{prop:extp} holds, and assume that 
$\sigma_{m,\vpi_{i}}(\pi)=\pi_{w_{1}\vpi_{i}} \ast 
\pi_{w_{2}\vpi_{i}} \ast \cdots \ast \pi_{w_{m}\vpi_{i}}$. 
Then we have $b:=x_{j_{1}}x_{j_{2}} \cdots x_{j_{k}}
u_{\vpi_{i}} \ne 0$, and 
$\sigma_{m,\vpi_{i}}(b)=u_{w_{1}\vpi_{i}} \otimes 
u_{w_{2}\vpi_{i}} \otimes \cdots \otimes u_{w_{m}\vpi_{i}}$. 
\end{lem}

%
\section{Main Results.}
\label{sec:main}

%
\subsection{Isomorphism theorem.}
\label{subsec:isom}

From now on, $\Fg$ is always an affine Lie algebra. 
We can carry out the proof of our isomorphism theorem, 
following the general line of that for \cite[Theorem 4.1]{Kas5}. 

\begin{thm} \label{thm:isom}
There exists a unique isomorphism $\Phi_{\vpi_{i}}:
\CB(\vpi_{i}) \stackrel{\sim}{\rightarrow} 
\BB_{0}(\vpi_{i})$ of crystals such that $\Phi_{\vpi_{i}}
(u_{\vpi_{i}})=\pi_{\vpi_{i}}$. 
\end{thm}

\begin{proof}
It suffices to prove that
for $j_{1},\,j_{2},\,\dots,\,j_{p} \in I$ and 
$k_{1},\,k_{2},\,\dots,\,k_{q} \in I$, 

\vspace{1.5mm}

\noindent (1) \ 
$x_{j_{1}}x_{j_{2}} \cdots x_{j_{p}}u_{\vpi_{i}}=
x_{k_{1}}x_{k_{2}} \cdots x_{k_{q}}u_{\vpi_{i}} \, 
\Leftrightarrow \, 
x_{j_{1}}x_{j_{2}} \cdots x_{j_{p}}\pi_{\vpi_{i}}=
x_{k_{1}}x_{k_{2}} \cdots x_{k_{q}}\pi_{\vpi_{i}}$, 

\vspace{1.5mm}

\noindent (2) \ 
$x_{j_{1}}x_{j_{2}} \cdots x_{j_{p}}u_{\vpi_{i}}=0 \, 
\Leftrightarrow \, 
x_{j_{1}}x_{j_{2}} \cdots x_{j_{p}}\pi_{\vpi_{i}} =0$. 

\vspace{1.5mm}

Part (2) has already been proved in Lemma~\ref{lem:extbp}.
Let us show the direction $(\Rightarrow)$ of 
Part (1). Take $m \in \BZ_{> 0}$ such that the assertion 
of Proposition~\ref{prop:extb} holds for 
both $b_{1}:=x_{j_{1}}x_{j_{2}} \cdots x_{j_{p}}u_{\vpi_{i}}$ 
and $b_{2}:=x_{k_{1}}x_{k_{2}} \cdots x_{k_{q}}u_{\vpi_{i}}$: 
\begin{align*}
& 
\sigma_{m,\vpi_{i}}(b_{1})=u_{w_{1}\vpi_{i}} \otimes 
u_{w_{2}\vpi_{i}} \otimes \cdots \otimes u_{w_{m}\vpi_{i}}, \\
& 
\sigma_{m,\vpi_{i}}(b_{2})=u_{w^{\prime}_{1}\vpi_{i}} \otimes 
u_{w^{\prime}_{2}\vpi_{i}} \otimes \cdots \otimes 
u_{w^{\prime}_{m}\vpi_{i}}.
\end{align*}
Since $b_{1}=b_{2}$, we get $u_{w_{l}\vpi_{i}}
=u_{w_{l}^{\prime}\vpi_{i}}$, and hence 
$w_{l}\vpi_{i}=w_{l}^{\prime}\vpi_{i}$ 
for all $l=1,\,2,\,\dots,\,m$. 
By Lemma~\ref{lem:extbp} (1), 
we see that 
\begin{align*}
& 
\sigma_{m,\vpi_{i}}(\pi_{1})=\pi_{w_{1}\vpi_{i}} \ast 
\pi_{w_{2}\vpi_{i}} \ast \cdots \ast \pi_{w_{m}\vpi_{i}}, \\
& 
\sigma_{m,\vpi_{i}}(\pi_{2})=\pi_{w^{\prime}_{1}\vpi_{i}} \ast 
\pi_{w^{\prime}_{2}\vpi_{i}} \ast \cdots \ast 
\pi_{w^{\prime}_{m}\vpi_{i}},
\end{align*}
where $\pi_{1}:=x_{j_{1}}x_{j_{2}} \cdots x_{j_{p}}\pi_{\vpi_{i}}$
and $\pi_{2}:=x_{k_{1}}x_{k_{2}} \cdots x_{k_{q}}\pi_{\vpi_{i}}$. 
Since $w_{l}\vpi_{i}=w_{l}^{\prime}\vpi_{i}$ and 
$\pi_{w\vpi_{i}}(t)=t(w\vpi_{i})$ for all $w \in W$, we get 
$\sigma_{m,\vpi_{i}}(\pi_{1})=\sigma_{m,\vpi_{i}}(\pi_{2})$. 
Since $\sigma_{m,\vpi_{i}}$ is injective, we conclude that 
$\pi_{1}=\pi_{2}$. 

We show the reverse direction $(\Leftarrow)$ of Part (1). 
Take $m \in \BZ_{> 0}$ such that the assertion of 
Proposition~\ref{prop:extp} holds for 
both $\pi_{1}:=x_{j_{1}}x_{j_{2}} \cdots x_{j_{p}}\pi_{\vpi_{i}}$ 
and $\pi_{2}:=x_{k_{1}}x_{k_{2}} \cdots x_{k_{q}}\pi_{\vpi_{i}}$: 
\begin{align*}
& 
\sigma_{m,\vpi_{i}}(\pi_{1})=\pi_{w_{1}\vpi_{i}} \ast 
\pi_{w_{2}\vpi_{i}} \ast \cdots \ast \pi_{w_{m}\vpi_{i}}, \\
& 
\sigma_{m,\vpi_{i}}(\pi_{2})=\pi_{w^{\prime}_{1}\vpi_{i}} \ast 
\pi_{w^{\prime}_{2}\vpi_{i}} \ast \cdots \ast 
\pi_{w^{\prime}_{m}\vpi_{i}}.
\end{align*}
Since $\pi_{1}=\pi_{2}$, and hence 
$\sigma_{m,\vpi_{i}}(\pi_{1})=
\sigma_{m,\vpi_{i}}(\pi_{2})$ in $\BP$, the two paths 
$\pi_{w_{1}\vpi_{i}} \ast 
\pi_{w_{2}\vpi_{i}} \ast \cdots \ast \pi_{w_{m}\vpi_{i}}$ 
and $\pi_{w^{\prime}_{1}\vpi_{i}} \ast 
\pi_{w^{\prime}_{2}\vpi_{i}} \ast \cdots \ast 
\pi_{w^{\prime}_{m}\vpi_{i}}$ are identical 
modulo reparametrization. Hence we can deduce that 
$w_{l}\vpi_{i}=w_{l}^{\prime}\vpi_{i}$ 
for all $l=1,\,2,\,\dots,\,m$ from the fact that 
if $a\vpi_{j} \in W\vpi_{i}$ for some $a \in \BQ_{\ge 0}$ and 
$i,\,j \in I_{0}$, then $i=j$ and $a=1$. 
By Lemma~\ref{lem:extbp} (2), we have 
\begin{align*}
& 
\sigma_{m,\vpi_{i}}(b_{1})=u_{w_{1}\vpi_{i}} \otimes 
u_{w_{2}\vpi_{i}} \otimes \cdots \otimes u_{w_{m}\vpi_{i}}, \\
& 
\sigma_{m,\vpi_{i}}(b_{2})=u_{w^{\prime}_{1}\vpi_{i}} \otimes 
u_{w^{\prime}_{2}\vpi_{i}} \otimes \cdots \otimes 
u_{w^{\prime}_{m}\vpi_{i}}.
\end{align*}
Since $w_{l}\vpi_{i}=w_{l}^{\prime}\vpi_{i}$ 
for all $l=1,\,2,\,\dots,\,m$, 
it follows from \cite[Proposition 5.8 (i)]{Kas5} that 
$u_{w_{l}\vpi_{i}}=u_{w_{l}^{\prime}\vpi_{i}}$ 
for all $l=1,\,2,\,\dots,\,m$. Therefore we have 
$\sigma_{m,\vpi_{i}}(b_{1})=
\sigma_{m,\vpi_{i}}(b_{2})$. 
Since $\sigma_{m,\vpi_{i}}$ is injective, 
we conclude that $b_{1}=b_{2}$. 
\end{proof}

In general, an isomorphism of crystals between 
$\CB(\lambda)$ and $\BB_{0}(\lambda)$ does not exist, 
even if $\CB(\lambda)$ is connected. For example, 
let $\Fg$ be of type $A_{2}^{(1)}$, and 
$\lambda=\vpi_{1}+\vpi_{2}$ (we know from 
\cite[Proposition 5.4]{Kas5} that $\CB(\lambda)$ is 
connected). If $\CB(\lambda) \cong \BB_{0}(\lambda)$ as crystals, 
then we would have $wu_{\lambda}=w^{\prime}u_{\lambda}$ 
in $\CB(\lambda)$ for every $w,\,w^{\prime} \in W$ with 
$w\lambda=w^{\prime}\lambda$, 
but we have an example of $w,\,w^{\prime} \in W$ 
such that $wu_{\lambda} \ne w^{\prime}u_{\lambda}$ 
in $\CB(\lambda)$ and $w\lambda=w^{\prime}\lambda$ 
(see \cite[Remark 5.10]{Kas5}). 
On the other hand, we obtain the following corollary by combining 
Theorems~\ref{thm:emb} and \ref{thm:isom}, because 
we can deduce from Theorem~\ref{thm:isom} that the concatenation 
$\Cat_{i \in I_{0}}\BB_{0}(\vpi_{i})^{\ast m_{i}}$
of $\BB_{0}(\vpi_{i})$'s is isomorphic to the tensor 
product $\bigotimes_{i \in I_{0}}\CB(\vpi_{i})^{\otimes m_{i}}$ 
of $\CB(\vpi_{i})$'s for all $m_{i} \in \BZ_{\ge 0},\, i \in I$.

\begin{cor}
Let $\lambda=\sum_{i \in I_{0}} m_{i}\vpi_{i}$, and 
set $\pi_{\ud{\lambda}}:=\Cat_{i \in I_{0}} 
\pi_{\vpi_{i}}^{\ast m_{i}} \in 
\Cat_{i \in I_{0}}\BB_{0}(\vpi_{i})^{\ast m_{i}}$.
Let $\BB_{0}(\ud{\lambda})$ be the connected component 
of $\Cat_{i \in I_{0}}\BB_{0}(\vpi_{i})^{\ast m_{i}}$ 
containing $\pi_{\ud{\lambda}}$, i.e., 
that of $\Cat_{i \in I_{0}}\BB(\vpi_{i})^{\ast m_{i}}$. 
Then there exists a 
unique isomorphism $\Phi_{\ud{\lambda}}:
\CB_{0}(\lambda) \stackrel{\sim}{\rightarrow} 
\BB_{0}(\ud{\lambda})$ of crystals 
that maps $u_{\lambda}$ to $\pi_{\ud{\lambda}}$. 
\end{cor}

\begin{rem}
In \cite{G}, Greenstein proved that 
if $\Fg$ is of type $A_{\ell}^{(1)}$, then the 
connected component $\BB_{0}(m\vpi_{i}+n\delta)$ is a 
path model for a certain bounded module $L(\ell,m,n)$. 
He also showed a decomposition rule for tensor products, 
which seems to be closely related to Theorem~\ref{thm:lr} below. 
\end{rem}

%
\subsection{Branching rule.}
\label{subsec:branch}

\begin{lem} \label{lem:bound}
For any $\pi \in \BB(\vpi_{i})$, we have 
$(\pi(1),\pi(1)) \le (\vpi_{i},\vpi_{i})$.
\end{lem}

\begin{proof}
Let $\pi=(\nu_{1},\,\nu_{2},\,\dots,\,\nu_{s} \, ; \, 
a_{0},\,a_{1},\,\dots,\,a_{s})$ with $\nu_{j} \in W\vpi_{i}$ 
and $a_{j} \in [0,1]$ be a Lakshmibai-Seshadri path of shape 
$\vpi_{i}$ (cf. \cite[\S4]{L2}). By the definition of 
a Lakshmibai--Seshadri path, we see that 
$\pi(1)=\sum_{j=1}^{s}(a_{j}-a_{j-1})\nu_{j}$. 
Hence we have 
\begin{align*}
(\pi(1),\pi(1)) & 
=\sum_{j=1}^{s}(a_{j}-a_{j-1})^{2}(\nu_{j},\nu_{j})+
2\sum_{1 \le k < l \le s}
(a_{k}-a_{k-1})(a_{l}-a_{l-1})(\nu_{k},\nu_{l}) \\
& =\sum_{j=1}^{s}(a_{j}-a_{j-1})^{2}(\vpi_{i},\vpi_{i})+
2\sum_{1 \le k < l \le s}
(a_{k}-a_{k-1})(a_{l}-a_{l-1})
(\vpi_{i},w_{kl}\vpi_{i})
\end{align*}
for some $w_{kl} \in W$. By \cite[Proposition 6.3]{Kac}, 
we deduce that $w_{kl}\vpi_{i}=\vpi_{i}-\beta_{kl}+n_{kl}\delta$ for 
some $\beta_{kl} \in \sum_{i \in I_{0}}\BZ_{\ge 0} \alpha_{i}$ 
and $n_{kl} \in \BZ$. Therefore we have 
(note that $\vpi_{i}$ is of level $0$)
\begin{align*}
(\pi(1),\pi(1)) & =\sum_{j=1}^{s}(a_{j}-a_{j-1})^{2}
(\vpi_{i},\vpi_{i}) \\[-1.0mm]
& \hspace*{20mm} 
+2\sum_{1 \le k < l \le s}
(a_{k}-a_{k-1})(a_{l}-a_{l-1})
(\vpi_{i},\vpi_{i}-\beta_{kl}+n_{kl}\delta) \\
&=\sum_{j=1}^{s}(a_{j}-a_{j-1})^{2}(\vpi_{i},\vpi_{i})+
2\sum_{1 \le k < l \le s}
(a_{k}-a_{k-1})(a_{l}-a_{l-1})
(\vpi_{i},\vpi_{i}) \\[-1.0mm]
& \hspace*{20mm}
-2\sum_{1 \le k < l \le s}
(a_{k}-a_{k-1})(a_{l}-a_{l-1})
(\vpi_{i},\beta_{kl}) \\
& = \biggl\{\sum_{j=1}^{s}(a_{j}-a_{j-1})\biggr\}^{2}
(\vpi_{i},\vpi_{i})-
2\sum_{1 \le k < l \le s}(a_{k}-a_{k-1})(a_{l}-a_{l-1})
(\vpi_{i},\beta_{kl}) \\
& =(\vpi_{i},\vpi_{i})-
2\sum_{1 \le k < l \le s}
(a_{k}-a_{k-1})(a_{l}-a_{l-1})
(\vpi_{i},\beta_{kl})
\end{align*}
Since $(\vpi_{i},\beta_{kl}) \ge 0$ for all $1 \le k < l \le s$, 
we deduce that $(\pi(1),\pi(1)) \le (\vpi_{i},\vpi_{i})$, as desired.
\end{proof}

Let $S$ be a proper subset of $I$, i.e., 
$S \subsetneq I$. Denote by $\Fg_{S}$ the Levi subalgebra of $\Fg$
corresponding to $S$. Note that a crystal $\BB$ for $\Fg$ 
can be regarded as a crystal for $\Fg_{S}$ by restriction.

\begin{thm} \label{thm:branch01}
As crystals for $\Fg_{S}$, $\BB(\vpi_{i})$ and $\BB_{0}(\vpi_{i})$ 
decompose as follows {\rm:}
\begin{equation} \label{eq:branch01}
\BB(\vpi_{i}) \cong \bigsqcup_{
\begin{subarray}{c}
\pi \in \BB(\vpi_{i}) \\
\pi\text{\rm: $\Fg_{S}$-dominant}
\end{subarray}}
\BB_{S}(\pi(1)), \qquad 
\BB_{0}(\vpi_{i}) \cong \bigsqcup_{
\begin{subarray}{c}
\pi \in \BB_{0}(\vpi_{i}) \\
\pi\text{\rm: $\Fg_{S}$-dominant}
\end{subarray}}
\BB_{S}(\pi(1)).
\end{equation}
where $\BB_{S}(\lambda)$ is the set of Lakshmibai--Seshadri paths
of shape $\lambda$ for $\Fg_{S}$, and a path $\pi$
is said to be $\Fg_{S}$-dominant if $(\pi(t))(\alpha^{\vee}_{i}) \ge 0$
for all $t \in [0,1]$ and $i \in S$.
\end{thm}

\begin{proof}
We will show the first equality in \eqref{eq:branch01}, 
since the second one can be shown in the same way.
As in \cite[\S9.3]{Kas1}, we deduce, using by Lemma~\ref{lem:bound}, 
that each connected component of $\BB(\vpi_{i})$ 
(as a crystal for $\Fg_{S}$) contains an extremal 
weight element $\pi^{\prime}$ with respect to 
$W_{S}:=\langle r_{j} \mid j \in S\rangle$. 
Because $\Fg_{S}$ is a finite-dimensional 
reductive Lie algebra, there exists 
$w \in W_{S}$ such that $((w\pi)(1))(\alpha^{\vee}_{j}) \ge 0$
for all $j \in S$. Put $\pi:=w\pi^{\prime}$ for this $w \in W_{S}$. 
Since $\pi$ is also extremal, we deduce from 
\cite[Lemma 4.5 d)]{L2} that $(\pi(t))(\alpha^{\vee}_{j}) \ge 0$
for all $t \in [0,1]$ and $j \in S$, i.e., 
$\pi$ is $\Fg_{S}$-dominant. We see from \cite[Theorem 7.1]{L2} that 
the connected component containing $\pi$ is isomorphic to 
$\BB_{S}(\pi(1))$, thereby completing the proof of the theorem. 
\end{proof}

\begin{thm} \label{thm:branch02}
{\rm (1)} \, 
The extremal weight module $V(\vpi_{i})$ of extremal weight $\vpi_{i}$ 
is completely reducible as a module over the quantized universal 
enveloping algebra $U_{q}(\Fg_{S})$ of the Levi subalgebra $\Fg_{S}$. 

\vspace{0.5mm}

\noindent
{\rm (2)} \, 
The decomposition of $V(\vpi_{i})$ as a $U_{q}(\Fg_{S})$-module 
is given by {\rm:}
\begin{equation}
V(\vpi_{i}) \cong \bigoplus_{
\begin{subarray}{c}
\pi \in \BB_{0}(\vpi_{i}) \\
\pi\text{\rm: $\Fg_{S}$-dominant}
\end{subarray}}
V_{S}(\pi(1)), 
\end{equation}
where $V_{S}(\lambda)$ is the integrable highest weight module 
of highest weight $\lambda$ over $U_{q}(\Fg_{S})$. 
\end{thm}

\begin{proof}
(1) We first prove that $U:=U_{q}(\Fg_{S})u$ is finite-dimensional
for each weight vector $u \in V(\vpi_{i})$. To prove this, 
it suffices to show that the weight system $\Wt(U)$ of $U$ 
is a finite set, since each weight space 
of $V(\vpi_{i})$ is finite-dimensional
(see \cite[Proposition 5.16 (iii)]{Kas5}). Remark that 
if $\mu,\,\nu \in P$ are weights of $U$, 
then $\mu,\,\nu \in \Fh^{\ast}_{0}$, and 
$\mu-\nu \in Q_{S}:=\sum_{i \in S}\BZ\alpha_{i}$. 
Hence the canonical map
$\cl:\Fh^{\ast}_{0} \twoheadrightarrow 
\Fh^{\ast}_{0}/\BQ\delta$ 
is injective on $\Wt(U)$, since $k\delta \not\in Q_{S}$ 
for any $k \in \BZ \setminus \{0\}$. 
Since $\Wt(U)$ is contained in the 
weight system $\Wt(V(\vpi_{i}))$ of $V(\vpi_{i})$, 
it follows from Theorem~\ref{thm:isom} and 
Lemma~\ref{lem:bound} that
\begin{align*}
\cl(\Wt(U)) & \subset \cl\bigl(\Wt(V(\vpi_{i}))\bigr)=
\cl\bigl(\{\pi(1) \mid \pi \in \BB_{0}(\vpi_{i})\}\bigr) 
\quad \text{by Theorem~\ref{thm:isom}} \\
& \subset \bigl\{\mu^{\prime} \in 
\Fh^{\ast}_{0}/\BQ\delta \mid 
(\mu^{\prime},\mu^{\prime}) \le 
(\cl(\vpi_{i}),\cl(\vpi_{i}))\bigr\} \quad 
\text{by Lemma~\ref{lem:bound}}. 
\end{align*}
Because the bilinear form $(\cdot\,,\,\cdot)$ on 
$\Fh^{\ast}_{0}/\BQ\delta$ is positive-definite, 
the set $\cl(\Wt(U))$ is discrete and contained 
in a compact set with respect to the metric topology
on $\BR \otimes_{\BQ} (\Fh^{\ast}_{0}/\BQ\delta)$ 
defined by $(\cdot\,,\,\cdot)$. Therefore, we see that 
$\cl(\Wt(U))$ is a finite set, and hence so is 
$\Wt(U)$. Thus, we conclude that 
$U=U_{q}(\Fg_{S})u$ is finite-dimensional. 

Since $q$ is assumed not to be a root of unity, the finite-dimensional 
$U_{q}(\Fg_{S})$-module $U_{q}(\Fg_{S})u$ is 
completely reducible for each weight vector $u \in V(\vpi_{i})$. 
Because $V(\vpi_{i})$ is a sum of all such $U_{q}(\Fg_{S})u$'s, 
we deduce that $V(\vpi_{i})$ is also completely reducible. 

(2) Because each weight space of 
$V(\vpi_{i})$ is finite-dimensional, we can define 
the formal character $\ch V(\vpi_{i})$ of $V(\vpi_{i})$. 
By Theorem~\ref{thm:branch01}, we have 
\begin{equation*}
\ch V(\vpi_{i})=\sum_{
\begin{subarray}{c}
\pi \in \BB_{0}(\vpi_{i}) \\
\pi\text{\rm: $\Fg_{S}$-dominant}
\end{subarray}}
\ch V_{S}(\pi(1)).
\end{equation*}
Therefore, in order to prove Part (2), 
we need only show that this is the unique way of writing 
$\ch V(\vpi_{i})$ as a sum of the characters 
of integrable highest weight 
$\Fg_{S}$-modules. Assume that
\begin{equation*}
\ch V(\vpi_{i})=\sum_{\lambda \in P}c_{\lambda}
\ch V_{S}(\lambda) \quad \text{and} \quad
\ch V(\vpi_{i})=\sum_{\lambda \in P}c^{\prime}_{\lambda}
\ch V_{S}(\lambda)
\end{equation*}
with $c_{\lambda},c^{\prime}_{\lambda} \in \BZ$ for $\lambda \in P$. 
Then we have $\sum_{\lambda \in P}
(c_{\lambda}-c^{\prime}_{\lambda})
\ch V_{S}(\lambda)=0$. Suppose that there exists $\lambda \in P$ 
such that $c_{\lambda}-c^{\prime}_{\lambda} \ne 0$, and 
set $X:=\bigl\{\lambda \in P \mid c_{\lambda}-c^{\prime}_{\lambda}
\ne 0\bigr\} (\ne \emptyset)$. Note that $X$ is contained in 
the weight system $\Wt(V(\vpi_{i}))$ of $V(\vpi_{i})$. 
As in the proof of Part (1), we deduce that
\begin{equation*}
\cl\bigl(\Wt(V(\vpi_{i}))\bigr) 
\subset \bigl\{\mu^{\prime} \in 
\Fh^{\ast}_{0}/\BQ\delta \mid 
(\mu^{\prime},\mu^{\prime}) \le 
(\cl(\vpi_{i}),\cl(\vpi_{i}))\bigr\},
\end{equation*}
and hence $\Wt(V(\vpi_{i}))$ 
modulo $\BZ\delta$ is a finite set.

Now, we define a partial order $\ge_{S}$ on $P$ as follows:
\begin{equation*}
\mu \ge_{S} \nu \quad
\text{for \ } \mu,\nu \in P 
\quad \Longleftrightarrow \quad \mu-\nu 
\in (Q_{S})_{+}:=\sum_{i \in S}\BZ_{\ge 0}\alpha_{i}.
\end{equation*}
Let us show that the set $X$ has a maximal element 
with respect to this order $\ge_{S}$. Let $\mu \in X$. 
Then $\Wt(V(\vpi_{i})) \cap (\mu+Q_{S})$ is a 
finite set. Indeed, if this is not a finite set, then 
there exist elements $\nu,\,\nu^{\prime}$ of it 
such that $\nu-\nu^{\prime}=k\delta$ with $k \in \BZ
\setminus \{0\}$, since $\Wt(V(\vpi_{i}))$ 
modulo $\BZ\delta$ is a finite set. However, 
since $\nu-\nu^{\prime} \in Q_{S}$ and $k\delta 
\not\in Q_{S}$ for any $k \in \BZ \setminus \{0\}$, 
this is a contradiction. Therefore, we see that
$X \cap (\mu+(Q_{S})_{+})$ is also a finite set, and hence 
that $X$ has a maximal element of the form $\mu+\beta$ 
for some $\beta \in (Q_{S})_{+}$. 

Let $\nu \in X$ be a maximal element 
with respect to this order $\ge_{S}$. 
We can easily see that 
the coefficient of $e(\nu)$ in 
$\sum_{\lambda \in P}
(c_{\lambda}-c^{\prime}_{\lambda})
\ch V_{S}(\lambda)$ is equal to $c_{\nu}-c^{\prime}_{\nu}$. 
Since $\nu \in X$, we have $c_{\nu}-c^{\prime}_{\nu} \ne 0$, 
which contradicts $\sum_{\lambda}
(c_{\lambda}-c^{\prime}_{\lambda})
\ch V_{S}(\lambda)=0$. This completes the 
proof of the theorem. 
\end{proof}

%
\subsection{Decomposition rule for tensor products.}
\label{subsec:lr}

In this subsection, we assume that $\vpi_{i}$ is minuscule, i.e., 
$\vpi_{i}(\alpha^{\vee}) \in \bigl\{\pm1,\,0\bigr\}$ for every dual 
real root $\alpha^{\vee}$ of $\Fg$. For example, if $\Fg$ is of type 
$A_{n}^{(1)}$, then all $\vpi_{i}$'s are minuscule. 

\begin{rem} \label{rem:minu}
If $\vpi_{i}$ is minuscule, then, 
for any $\mu,\,\nu \in W\vpi_{i}$
and rational number $0 < a < 1$, there does not 
exist an $a$-chain for $(\mu,\nu)$. Hence it follows from the 
definition of a Lakshmibai--Seshadri path that 
$\BB(\vpi_{i})=\bigl\{\pi_{w\vpi_{i}} \mid w \in W\bigr\}$. 
Since $w\pi_{\vpi_{i}}=\pi_{w\vpi_{i}}$, we see that 
$\BB(\vpi_{i})$ is connected, and hence 
$\BB(\vpi_{i})=\BB_{0}(\vpi_{i})$. 
\end{rem}

\begin{thm} \label{thm:lr}
Let $\lambda$ be a dominant integral weight which is not 
a multiple of the null root $\delta$ of $\Fg$. Then 
the concatenation $\BB(\lambda) \ast \BB(\vpi_{i})$ 
decomposes as follows\,{\rm:}
\begin{equation}
\BB(\lambda) \ast \BB(\vpi_{i})
 \cong {\displaystyle\bigsqcup_{
\begin{subarray}{c}
\pi \in \BB(\vpi_{i}) \\
\pi\text{\rm: $\lambda$-dominant}
\end{subarray}}}
\BB(\lambda+\pi(1)),
\end{equation}
where $\pi \in \BB(\vpi_{i})$ is said to be 
$\lambda$-dominant if $(\lambda+\pi(t))(\alpha^{\vee}_{i}) \ge 0$
for all $t \in [0,1]$ and $i \in I$.
\end{thm}

\begin{proof}
We will prove that each connected component contains 
a (unique) path of the form $\pi_{\lambda} \ast \pi$
for a $\lambda$-dominant path $\pi \in \BB(\vpi_{i})$. 
Then the assertion of the theorem 
follows from \cite[Theorem 7.1]{L2}. 

Let $\pi_{1} \ast \pi_{2} \in 
\BB(\lambda) \ast \BB(\vpi_{i})$. It can easily be seen 
that $e_{i_{1}}e_{i_{2}} \cdots e_{i_{k}} 
(\pi_{1} \ast \pi_{2})=\pi_{\lambda} \ast \pi_{2}^{\prime}$ 
for some $i_{1},\,i_{2},\,\dots,\,i_{k} \in I$, where 
$\pi_{2}^{\prime} \in \BB(\vpi_{i})$ (cf. \cite[\S5.6]{G}).
Set $S:=\bigl\{i \in I \mid 
\lambda(\alpha^{\vee}_{i})=0\bigr\}$ (note that 
$S \subsetneq I$, since $\lambda$ is not a 
multiple of $\delta$), and 
let $\BB$ be the set of paths of the form 
$e_{j_{1}}e_{j_{2}} \cdots e_{j_{l}}
(\pi_{\lambda} \ast \pi_{2}^{\prime})$
for some $j_{1},\,j_{2},\,\dots,\,j_{l} \in S$. 
Remark that if $e_{j_{1}}e_{j_{2}} \cdots e_{j_{l}}
(\pi_{\lambda} \ast \pi_{2}^{\prime}) \ne 0$, 
then $e_{j_{1}}e_{j_{2}} \cdots e_{j_{l}}
(\pi_{\lambda} \ast \pi_{2}^{\prime})=\pi_{\lambda} \ast 
(e_{j_{1}}e_{j_{2}} \cdots e_{j_{l}}\pi_{2}^{\prime})$. 
As in the proof of Part (2) of Theorem~\ref{thm:branch02}, 
we deduce that 
\begin{equation*}
\bigl\{\pi(1) \mid \pi \in \BB(\vpi_{i})\bigr\} 
\cap \bigl(\pi_{2}^{\prime}(1)+(Q_{S})_{+}\bigr)
=\Wt(V(\vpi_{i})) \cap \bigl(\pi_{2}^{\prime}(1)+(Q_{S})_{+}\bigr)
\end{equation*}
is a finite set. Hence we have $\pi_{\lambda} \ast 
\pi_{2}^{\prime\prime} \in \BB$ for some $\pi_{2}^{\prime\prime} 
\in \BB(\vpi_{i})$ such that $e_{j}
(\pi_{\lambda} \ast \pi_{2}^{\prime\prime})=0$, $j \in S$. Because $\vpi_{i}$ is minuscule and 
$\pi_{2}^{\prime\prime}=\pi_{w\vpi_{i}}$ for some 
$w \in W$ (cf. Remark~\ref{rem:minu}), we see that 
$e_{j}(\pi_{\lambda} \ast \pi_{2}^{\prime\prime})=0$ 
for all $j \in I \setminus S$. Therefore,
we conclude that $\pi_{2}^{\prime\prime} \in \BB(\vpi_{i})$ is 
$\lambda$-dominant. Thus we have completed the proof 
of the theorem. 
\end{proof}

\begin{rem}
Unlike Theorems~\ref{thm:branch01} and \ref{thm:branch02}, this 
theorem does not necessarily imply the decomposition rule for tensor 
products of corresponding $U_{q}(\Fg)$-modules.
\end{rem}


{\small
\setlength{\baselineskip}{12pt}
\renewcommand{\refname}{References}

}

\end{document}